\documentclass[a4paper,12pt,reqno,twoside]{amsart}

\usepackage[utf8]{inputenc} 		% UTF8 encoding
\usepackage[T1]{fontenc} 
\usepackage{mathpazo}

\usepackage{amssymb}
\usepackage{amsmath}
\usepackage{graphicx}
\usepackage[colorlinks,citecolor=blue,urlcolor=blue]{hyperref}
\usepackage{cite}
\usepackage[hmargin=2.5cm,vmargin=2.5cm]{geometry}
\usepackage{mathrsfs} % Script fonts

\newenvironment{myproof}{\paragraph{Proof:}}{\hfill$\blacksquare$}

\newcommand{\be}{\begin{eqnarray}}
\newcommand{\ee}{\end{eqnarray}}

\newcommand{\beq}{\begin{equation}}
\newcommand{\eeq}{\end{equation}}
\newcommand{\beqn}{\begin{equation*}}
\newcommand{\eeqn}{\end{equation*}}

\newcommand{\round}[1]{\lfloor#1\rfloor}

\newtheorem{thm}{Theorem}[section]

\newtheorem{cor}[thm]{Corollary}
\newtheorem{lem}[thm]{Lemma}
\newtheorem{defn}[thm]{Definition}
\newtheorem{fact}[thm]{Fact}

\theoremstyle{remark}
\newtheorem{remark}[thm]{Remark}

\newcommand\cC{{\mathcal C}}

\newcommand\cL{{\mathcal L}}
\newcommand\cM{{\mathcal M}}

\newcommand\cP{{\mathcal P}}

\newcommand\bE{{\mathbb E}}

\newcommand\bP{{\mathbb P}}

\newcommand\bR{{\mathbb R}}

\def\bfT{\mathbf{T}}

\renewcommand{\qed}{\hfill\blacksquare}

\begin{document}

\title[Intermittent quasistatic dynamical systems: convergence of fluctuations]{Intermittent quasistatic dynamical systems: \\weak convergence of fluctuations}

\author[Juho Lepp\"anen]{Juho Lepp\"anen}
\address[Juho Lepp\"anen]{
Department of Mathematics and Statistics, P.O.\ Box 68, Fin-00014 University of Helsinki, Finland.}
\email{juho.leppanen@helsinki.fi}
%\urladdr{http://www.math.helsinki.fi/mathphys/mikko.html}

\keywords{Quasistatic dynamical system, intermittency, Pomeau--Manneville map, weak convergence, martingale problem, functional central limit theorem}

\thanks{2010 {\it Mathematics Subject Classification.} 37C60; 37D25, 37A10, 60G44, 60H10} % Suggesting these. 

\maketitle

% Suggesting these. 
 
% http://www.ams.org/msc/msc2010.html
% 37A10  	One-parameter continuous families of measure-preserving transformations
% 37A15  	General groups of measure-preserving transformations [See mainly 22Fxx]
% 37A17  	Homogeneous flows [See also 22Fxx]
% 37A20  	Orbit equivalence, cocycles, ergodic equivalence relations
% 37A25  	Ergodicity, mixing, rates of mixing
% 37A30  	Ergodic theorems, spectral theory, Markov operators
% 37D20  	Uniformly hyperbolic systems (expanding, Anosov, Axiom A, etc.)
% 37D25  	Nonuniformly hyperbolic systems (Lyapunov exponents, Pesin theory, etc.)
% 37C40  	Smooth ergodic theory, invariant measures
% 37C60  	Nonautonomous smooth dynamical systems
% 60G44  Martingales with continuous parameter
% 60H10  	Stochastic ordinary differential equations
\begin{abstract} This paper is about statistical properties of quasistatic dynamical systems. These are a class of non-stationary systems that model situations where the dynamics change very slowly over time due to external influences. We focus on the case where the time-evolution is described by intermittent interval maps (Pomeau-Manneville maps) with time-dependent parameters. In a suitable range of parameters, we obtain a description of the statistical properties as a stochastic diffusion, by solving a well-posed martingale problem. The results extend those of a related recent study due to Dobbs and Stenlund, which concerned the case of quasistatic (uniformly) expanding systems.
\end{abstract}

\section{Introduction}

In this paper we continue the study, initiated in \cite{leppanen2016}, of intermittent quasistatic dynamical systems. These are non-uniformly expanding examples of the following non-stationary systems recently introduced by Dobbs and Stenlund \cite{dobbs2016}.

\begin{defn}[Discrete time QDS]\label{def:QDS}
Let~$(X,\mathscr{F})$ be a measurable space, $\cM$~a topological space whose elements are measurable self-maps $T:X\to X$, and~$\bfT$ a triangular array of the form
\beqn\label{eq:array}
\bfT = \{T_{n,k}\in\cM\ :\ 0\le k\le n, \ n\ge 1\} .
\eeqn
If there exists a piecewise continuous curve $\tau:[0,1]\to\cM$ such that\footnote{For any real number \(s \ge 0\), \(\round{s}\) denotes the integer part of \(s\).}
\begin{align}\label{gamma_speed}
\lim_{n\to\infty}T_{n,\round{nt}} = \tau_t
\end{align}
for all $t$,
we say that $(\bfT, \tau)$ is a \emph{quasistatic dynamical system (QDS)} with \emph{state space}~$X$ and \emph{system space}~$\cM$. 
\end{defn}

QDSs model situations where external influences force the observed system to transform very gradually over time; see \cite{stenlund2016,dobbs2016} for more discussions on their physical interpretation and significance. The compositions \(T_{n,k} \circ \cdots \circ T_{n,1}\) typically fail to be identically distributed. Such systems, which lack invariant measure, have gained interest during recent years due to advances in the research of time-dependent dynamical systems; see for example \cite{tanzi2016, conze2007,freitas2016,gupta2013,haydn2017,kawan2014, kawan2015, korepanov2016,mohapatra2014,nandori2012,nicol2016,
stenlund2014,ott2009, rousseau2016}.\\
\indent In a QDS $(\bfT, \tau)$ , the time-evolution of a point \(x \in X\) is described by the array \(\bfT\), separately on each level of the array: For each \(n \ge 1\), \(x_{n,k} = T_{n,k} \circ \cdots \circ T_{n,1} (x) \) is the state of the system after \(k \le n\) steps on the \(n\)th level. We are interested in the statistical properties of \( (x_{n,k})_{0 \le k \le n} \) as \(n \to \infty\). These depend on the limiting curve \(\tau\), which is approximated by the piecewise constant curve \(t \mapsto T_{n,\round{nt}}\) with increasing accuracy as \(n\) grows. In the analysis of a particular system, it might be necessary to specify the rate of convergence in \eqref{gamma_speed}. This is the case with the system studied in the present paper (defined below), along with those considered in \cite{dobbs2016,stenlund2016}.\\
\indent Given a measurable function \(f: X \to \bR\), and integers \(k,n\) with \(0 \le k \le n\), we denote \(f_{n,k} = f \circ T_{n,k} \circ \cdots \circ T_{n,1}\) and adopt the convention that \(f_{n,0} = f\). We define the functions \(S_n : X \times [0,1] \to \bR \) by
\begin{align*}
S_n(x,t) = \int^{nt}_0 f_{n,\round{s}}(x) \, ds, \hspace{0.5cm} n \ge 1.
\end{align*}
Note that for any \(x \in X\), the map \(t \to S_n(x,t)\) is a piecewise linear interpolation of the Birkhoff-type sum \(\sum_{k=0}^{\round{nt}} f_{n,k}(x) \) and thus belongs to the space of continuous functions \([0,1] \to \bR\), which we denote by \(C([0,1])\). Given an initial distribution \(\mu\) for \(x \in X\), we equip \(C([0,1])\) with the uniform norm topology, and view each map \(x \mapsto S_n(x,\cdot)\) as a random element with values in \(C([0,1])\).\\
\indent We now briefly discuss the results of \cite{dobbs2016, stenlund2016, leppanen2016}, and explain their connections to the current manuscript. In \cite{stenlund2016}, Stenlund showed an ergodic theorem for a general QDS, and applied the result to a particular model $(\bfT, \tau)$ whose systems space consists of strongly chaotic (uniformly) expanding maps on the circle \(X = S^1\), \(\tau\) is piecewise Hölder continuous, and the convergence rate in \eqref{gamma_speed} is sufficiently rapid (for details of the model, see Section 3 of \cite{stenlund2016}). The application yielded for any continuous function \(f : X \to [0,1]\), and almost every \(x \in X\) (with respect to Lebesgue measure),
\begin{align}\label{ergodic}
\lim_{n\to\infty} \sup_{t \in [0,1]} | n^{-1}S_n(x,t) - \int_0^t \hat{\mu}_{\tau_s} (f) \, ds | = 0,
\end{align} 
where \(\hat{\mu}_{\tau_s}\) denotes the unique SRB measure (equivalent to \(m\)) associated to \(\tau_s \in \cM\), and \(\hat{\mu}_{\tau_s}(f) = \int_0^1 f \, d\hat{\mu}_{\tau_s}\). An ergodic theorem for quasistatic billiards was also proved in \cite{stenlund2016}. \\
\indent In \cite{leppanen2016}, our primary aim was to investigate the ergodic properties of a QDS whose system space consists of non-uniformly expanding maps, and in particular extend the result \eqref{ergodic} to this setting. To this end, we introduced an intermittent version of the QDS, whose definition we next recall.\\
\indent For each \(\alpha \in [0,1)\), let \(T_{\alpha} : [0,1] \to [0,1]\) be the Pomeau-Manneville map from \cite{liverani1999} defined by
\begin{align*}
T_{\alpha }(x) = \begin{cases} x(1+ 2^{\alpha }x^{\alpha}) & \forall x \in [0, 1/2), \\
2x-1 & \forall x \in [1/2,1].
 \end{cases}
\end{align*}

\begin{defn}[Intermittent QDS]\label{def_int_qds}
Let~$X = [0,1]$ and $\cM = \{T_\alpha\,:\,0 \le\alpha < 1\}$ (equipped, say, with the uniform topology). Next, let
\beqn
\{\alpha_{n,k}\in [0,1)\ :\ 0\le k\le n, \ n\ge 1\}
\eeqn
be a triangular array of parameters and
\beqn
\gamma : [0,1] \to [0,1)
\eeqn
a piecewise continuous curve satisfying
\beqn
\lim_{n\to\infty} \alpha_{n,\round{nt}} = \gamma_t
\eeqn
for all $t$.
Finally, define $\tau_t = T_{\gamma_t}$ and
\beqn\label{eq:T_nk}
T_{n,k} = T_{\alpha_{n,k}}.
\eeqn
\end{defn}

For \(\alpha = 0\), \(T_{\alpha}\) is the angle-doubling map. For \(\alpha > 0\), \(T_{\alpha}\) is strongly chaotic due to expansion except near the neutral fixed point at the origin. The neutrality of this fixed point is determined by \(\alpha\): the larger the \(\alpha\), the longer it takes for points to escape a small neighborhood of the origin. It is well-known \cite{aimino2015,liverani1999} that every \(T_{\alpha}\) admits an invariant absolutely continuous probability measure \(\hat{\mu}_{\alpha}\), whose density belongs to the convex cone
\beqn
\begin{split}
\cC_*(\alpha) = \{f\in C((0,1])\cap L^1\,:\, & \text{$f\ge 0$, $f$ decreasing,} 
\\
& \text{$x^{\alpha+1}f$ increasing, $f(x)\le 2^{\alpha} (2 + \alpha) x^{-\alpha} m(f)$}\},
\end{split}
\eeqn
where \(m(f) = \int_0^1 f(x) \, dx \). We denote this density by \(\hat{h}_{\alpha}\). \\
\indent The main result of \cite{leppanen2016} showed that \eqref{ergodic} continues to hold for the above intermittent QDS, when the limiting curve \(\gamma\) is piecewise Hölder-continuous of order \(\eta \in (0,1]\), such that
\begin{align}\label{gamma_as_1_intro}
\overline{\gamma([0,1])} \subset [0,\beta_*]
\end{align}
holds for some \(\beta_* < \tfrac12\), and 
\begin{align}\label{gamma_as_2_intro}
\lim_{n\to\infty} n^{\eta}\sup_{t \in [0,1]} | \alpha_{n,\round{nt}} - \gamma_t | < \infty.
\end{align}
The former condition enables one to maintain uniform control in estimates involving correlation decay, while the latter condition reflects the regularity of \(\gamma\). Note that \eqref{gamma_as_2_intro} is always satisfied by the "equipartition" \(\alpha_{n,k} = \gamma_{kn^{-1}}\). If the assumption on \(\beta_*\) is relaxed to \(\beta_* < 1\), we managed to show in \cite{leppanen2016} that \eqref{ergodic} still holds if almost sure convergence is replaced by the weaker notion of convergence in probability (with respect to Lebesgue measure).\\
\indent Besides the aforementioned ergodic properties, other statistical properties of the intermittent QDS have not been widely studied. In this paper we consider distributional properties of the paths \(t \mapsto S_n(x,t)\) in the case where \eqref{gamma_as_1_intro} holds with $\beta_* < \tfrac12$, given an initial distribution \(\mu\) of \(x \in X\) with cone density. We obtain results which show that, for a wide class of centering measures \(\nu\), the fluctuations
\begin{align*}
\chi_n^{\nu}(x,t) &= n^{-\frac12}S_n(x,t) - \nu (n^{-\frac12}S_n(x,t))
\end{align*}
converge weakly to a stochastic diffusion process. These results extend those of \cite{dobbs2016} concerning the previously mentioned uniformly expanding QDS. We require \(f\) to be Lipschitz continuous, and the limiting curve \(\gamma\) to be Hölder-continuous such that \eqref{gamma_as_2_intro} along with \eqref{gamma_as_1_intro} holds. Additionally, we need to assume that the \(\mu\)-centered sequence \((\chi_n^{\mu})\) is tight. In the smaller parameter range \(\beta_* < \tfrac13\), we show that tightness holds for all Lipschitz continuous functions \(f\).

\subsection{Main result.} We work in the setting of Definition \ref{def_int_qds}:
\begin{align*}
\{ T_{n,k} \ :\ 0\le k\le n, \ n\ge 1\}
\end{align*}
is a fixed triangular array of Pomeau-Manneville maps \(T_{\alpha_{n,k}} = T_{n,k}\), and \(\gamma\) is a continuous curve \([0,1] \to [0,1)\) such that \(\lim_{n\to\infty} \alpha_{n,\round{nt}} = \gamma_t\). Recall that
\begin{align*}
S_n(x,t) = \int^{nt}_0 f_{n,\round{s}}(x) \, ds,
\end{align*}
where \(f_{n,k} = f \circ T_{n,k} \circ \cdots \circ T_{n,1}\). We fix some centering measure \(\nu\), denote
\begin{align*}
\bar{f} = f - \nu(f),
\end{align*}
\indent and for each integer \(n \ge 1\) define the fluctuation \(\chi_n^{\nu} : [0,1] \times [0,1] \to \bR\) by
\begin{align*}
\chi_n^{\nu}(x,t) &= n^{-\frac12}S_n(x,t) - n^{-\frac12}\nu(S_n(x,t)) \\
&= n^{\frac12} \int_0^t f_{n,\round{ns}}(x) \, ds - n^{\frac12} \int_0^t \nu(f_{n,\round{ns}}(x)) \, ds \\
&= n^{\frac12} \int_0^t \bar{f}_{n,\round{ns}}(x) \, ds.
\end{align*}
For brevity we usually hide the \(x\)-dependence here and denote \(\chi_n^{\nu}(t) = \chi_n^{\nu}(x,t)\). Given an initial probability measure \(\mu\), the map \(x \mapsto \chi_n^{\nu}(x,\cdot)\) is a random element with values in \(C([0,1])\), and we denote its distribution (with respect to \(\mu\)) by \(\bP^{\mu,\nu}_n\). If the centering measure \(\nu = \mu\), we denote \(\bP^{\mu,\mu}_n = \bP^{\mu}_n\). \\
\indent Recall that \(\hat{\mu}_{\alpha}\) denotes the invariant SRB measure associated to \(T_{\alpha}\). For all \(t \in [0,1]\), we set
\begin{align*}
\hat{f}_t = f - \hat{\mu}_{\gamma_t}(f),
\end{align*}
and 
\begin{align*}
\hat{\sigma}^2_t(f) = \lim_{m\to \infty} \hat{\mu}_{\gamma_t} \left[ \left( \frac{1}{\sqrt{m}} \sum_{k=0}^{m-1} \hat{f}_t \circ T_{\gamma_t}^k \right)^2 \right].
\end{align*}
In other words, \(\hat{\sigma}^2_t(f)\) is the limiting variance of \(\frac{1}{\sqrt{m}} \sum_{k=0}^{m-1} \hat{f}_t \circ T_{\gamma_t}^k \) with respect to the measure \(\hat{\mu}_{\gamma_t}\). \\
\indent We now come to the main result of this article:

\begin{thm}\label{main} Let \(f: \, [0,1] \to \bR\) be Lipschitz continuous, and let the initial measure \(\mu\)  be such that its density belongs to  \(\cC_*(\beta_*)\). Suppose that \(\gamma : [0,1] \to [0,1)\) is Hölder-continuous of order \(\eta \in (0,1]\), that \(\gamma([0,1]) \subset [0,\beta_*]\) for some \(\beta_* < \tfrac12\), and that 
\begin{align*}
\lim_{n\to\infty} n^{\eta}\sup_{t \in [0,1]} | \alpha_{n,\round{nt}} - \gamma_t | < \infty.
\end{align*}
Then, the variance \(\hat{\sigma}^2_t(f)\) is finite and satisfies the Green-Kubo formula
\begin{align*}
\hat{\sigma}^2_t(f) = \hat{\mu}_{\gamma_t}[\hat{f}^2_t] + 2 \sum_{k=1}^{\infty} \hat{\mu}_{\gamma_t}[\hat{f}_t \hat{f}_t \circ T_{\gamma_t}^k].
\end{align*}
If the sequence of measures \((\bP_n^{\mu})_{n\ge 1}\) is tight, then for any probability measure \(\nu\), whose density \(g=g_1 - g_2\) for some \(g_1,g_2 \in \cC_*(\beta_*)\), the sequence \((\bP_n^{\mu,\nu})_{n\ge 1}\) converges weakly to the law of the process
\begin{align}\label{diff}
\chi(t) = \int_0^t \hat{\sigma}_s(f) \, dW_s.
\end{align}
Here \(W\) is a standard Brownian motion, and the stochastic integral is to be understood in the sense of It\(\overline{o}\).
\end{thm}

\begin{remark} A couple of remarks are in order:
\begin{itemize}
\item[(1)]{If the density \(g \in \cC_*(\beta_*)\) or if \(g\) is Lipschitz continuous, then by Lemma 2.4 in \cite{leppanen2017} there exist \(g_1, g_2 \in \cC_*(\beta_*)\) such that \(g = g_1 - g_2\).}\smallskip
\item[(2)]{If \(f : [0,1] \to \bR^d\) is a vector-valued Lipschitz continuous function, Theorem \ref{main} continues to hold with the modifications that \(W\) is a \(d\)-dimensional standard Brownian motion, and \(\hat{\sigma}_t(f) \in \bR^{d \times d}\) is the square root of the covariance matrix  
\begin{align*}
\hat{\sigma}^2_t(f) = \lim_{m \to \infty} \mu_{\gamma_t} \left[ \left( \frac{1}{\sqrt{m}} \sum_{k=0}^{m-1} \hat{f}_t \circ T_{\gamma_t}^k \right) \otimes \left( \frac{1}{\sqrt{m}} \sum_{k=0}^{m-1} \hat{f}_t \circ T_{\gamma_t}^k \right) \right].
\end{align*}
Here \(v\otimes w\) is the \(d \times d\)-matrix with entries \((v\otimes w)_{\alpha\beta} = v_\alpha w_\beta\). To obtain this generalization,  it suffices to modify the proof of Theorem \ref{main} exactly as described in Section 9.1 of \cite{dobbs2016}.}
\end{itemize}
\end{remark}

A result similar to Theorem \ref{main} was established in \cite{dobbs2016} for a class of uniformly expanding QDSs. To prove the above theorem, we closely follow the approach of \cite{dobbs2016} and identify the limit process \(\chi\) by solving a well-posed martingale problem. We need tightness for the sequence \((\bP_n^{\mu})_{n \ge 0} \) to ensure the existence of a weakly convergent subsequence. If \(\beta_*\) is sufficiently small, this follows easily from correlation decay.

\begin{cor} Let \(f: \, [0,1] \to \bR\) be Lipschitz continuous, and let \(\gamma\) and \(\beta_*\) be as in Theorem \ref{main}. Assume also that \(\beta_* < \tfrac13\). Then, for any measure \(\mu\) whose density belongs to \(\cC_*(\beta_*)\), and for any measure \(\nu\) whose density \(g=g_1 - g_2\) for some \(g_1,g_2 \in \cC_*(\beta_*)\), the sequence of measures \((\bP_n^{\mu,\nu})_{n\ge 1}\) converges weakly to the law of the process
\begin{align*}
\chi(t) = \int_0^t \hat{\sigma}_s(f) \, dW_s.
\end{align*}
\end{cor}

\begin{myproof} By Theorem \ref{main}, it suffices to show that \((\bP_n^{\mu})_{n\ge 1}\)  is tight. This is verified in Lemma \ref{tight}. \end{myproof}

\indent Although on a general level, the proof of Theorem \ref{main} proceeds exactly as that of \cite{dobbs2016}, on the level of details the two proofs differ significantly. Often these differences are related to the fact that  exponential memory loss, which was a key ingredient in the proof of \cite{dobbs2016}, fails for the intermittent QDS. Instead we have polynomial memory loss (given by Fact \ref{aimino} below), and in the setting of Theorem \ref{main} the polynomial rate is determined by \(\beta_*\).   If tightness of \((\bP_n^{\mu})_{n\ge 1}\) is given, we manage to work in the parameter range \(\beta_* < \tfrac12\), albeit considerable subtlety is often required to deal with estimates involving correlation decay. Moreover, we are unable to prove that \((\bP_n^{\mu})_{n\ge 1}\) is tight for a large class of observables \(f\), if \(\tfrac13 \le \beta_* < \tfrac12\). To discuss these issues in more detail, we proceed to outline the proof of the theorem. 

\subsection{Outline of the proof of Theorem \ref{main}.} In the proof we extensively utilize the results of \cite{aimino2015,leppanen2017} concerning polynomial memory loss of time-dependent intermittent maps, and the results of \cite{leppanen2016} on perturbation of transfer operators. These are reviewed in Section \ref{prelim}. \\
\indent We start by denoting \(\xi_n = \chi^{\mu}_n\) and observe that it suffices to prove Theorem \ref{main} in the case where the centering measure \(\nu = \mu\), for an argument utilizing the Portmanteau theorem and polynomial memory loss then implies the more general result. If \((\bP_n^{\mu})_{n \ge 0} \) is tight, we know that it has a weakly convergent subsequence. By Kolmogorov-Chentsov criterion, tightness follows from the fourth moment bound
\begin{align}\label{tight_intro}
\mu[[\xi_n(t + \delta) - \xi_n(t)]^4] \lesssim \Vert f \Vert_{\text{Lip}}^4\delta^2.
\end{align}
In \cite{dobbs2016}, a bound of the form \eqref{tight_intro} was a direct consequence of exponential correlation decay. We invoke polynomial correlation decay (Fact \ref{multicor} below), which suffices to show \eqref{tight_intro} for all Lipschitz continuous \(f\), if \(\beta_* < \tfrac13\). If instead \(\tfrac13 \le \beta_* < \tfrac12\), we obtain a bound weaker than \eqref{tight_intro}, which still suffices for the rest of the proof, but in this case we have to assume that \((\bP_n^{\mu})_{n \ge 0} \) is tight to guarantee the existence of a weakly convergent subsequence. \\
\indent The existence of a weak limit \(\bP = \lim_{k} \bP_{n_k}^{\mu}\) along a subsequence enables us to pursue the path taken in \cite{dobbs2016}: We show that \(\bP\) solves the martingale problem corresponding to the expression of \(\chi\) in \eqref{diff}. By uniqueness of such solutions, the result of Theorem \ref{main} follows. \\
\indent A successful implementation of the above strategy requires control over the second moment \(\mu[[\xi_n(t) - \xi_n(s)]^2]\). Following \cite{dobbs2016}, we prove a representation
\begin{align}\label{intro_approx2}
\mu[[\xi_n(t+\delta) - \xi_n(t)]^2] =  \int_{t}^{t+\delta} \hat{\sigma}^2_s(f) \, ds + \delta o(1) +  o(n^{-\frac12}).
\end{align}
The proof of this result rests on the observation that \(\mu(\bar{f}_{n,\round{ns}}\bar{f}_{n,\round{nr}})\) is small outside a neighborhood \(A_n\) of the diagonal \(\{ (s,r) \: : \: t \le s = r \le t + \delta \} \), so that for large \(n\) we can approximate
\begin{align*}
\mu[[\xi_n(t+\delta) - \xi_n(t)]^2] \approx n \iint_{A_n} \mu(\bar{f}_{n,\round{ns}}\bar{f}_{n,\round{nr}}) \, dr \,ds.
\end{align*}
A result of \cite{leppanen2016} implies that for some \(\kappa, \theta \in (0,1)\),
\begin{align}\label{intro_approx5}
\mu(f_{n,\round{nr}}) = \hat{\mu}_{\gamma_r}(f) + O(n^{-\theta})
\end{align}
holds whenever \(nr \gtrsim  n^{\kappa} \). In order to make use of this we remove small (polynomially decaying) blocks from the lower left and upper right corners of \(A_n\), but here we need to be more careful than in \cite{dobbs2016} when choosing the size of these blocks: In the QDS of \cite{dobbs2016}, the weaker lower bound \(nr \gtrsim \log n \) sufficed for an estimate similar to \eqref{intro_approx5}. After these steps, we arrive at \eqref{intro_approx2} by using the perturbation estimates of  \cite{leppanen2016}, which allow us to approximate
\begin{align*}
\mu(f_{n,\round{ns}} f_{n,\round{nr}}) \approx \mu(f \circ T_{\gamma_s}^{\round{ns}} f \circ T_{\gamma_s}^{\round{nr}}) \hspace{0.5cm} \forall (s,r) \in A_n.
\end{align*}
 \indent Another auxiliary result instrumental in the proof of Theorem \ref{main} is the following decorrelation estimate for \(\xi_n\): whenever \(A: \bR \to \bR\) is a bounded function, and \(s \le t\), 
\begin{align}\label{intro_approx3}
 \mu[A(\xi_n(s))[\xi_n(t) - \xi_n(s)]^2] = \mu[A(\xi_n(s)]\mu[[\xi_n(t) - \xi_n(s)]^2] + o(1),
\end{align}
as $n \to \infty$. To prove the result, we  start by introducing (for each \(n\) and \(t\)) a canonical partition \(\cP_{n,t}\) of the unit interval, such that the map \(x \mapsto \xi_n(x,t)\) is almost constant on each partition element. This enables us to reduce \eqref{intro_approx3} to the following statement: For each partition element \(I \in \cP_{n,t}\), letting \(\mu_{I}\) denote the measure \(\mu\) conditioned on \(I\),
\begin{align}\label{intro_approx4}
\sum_{I \in \cP_{n,t}} \mu(I) \mu_{I}[[ \xi_n(t) - \xi_n(s)]^2] = \mu[[\xi_n(t) - \xi_n(s)]^2] + o(1),
\end{align}
as $n \to \infty$. In the uniformly expanding framework of \cite{dobbs2016}, the authors proved \eqref{intro_approx4} by observing that, if a suitably small block \( [s,s + n^{-p} ]^2\subset [s,t]^2\) is removed from the domain of integration, exponential loss of memory implies a uniformly small upper bound to
\begin{align}\label{intro_approx6}
\mu_I[[ \xi_n(t) - \xi_n(s)]^2] - \mu[[\xi_n(t) - \xi_n(s)]^2].
\end{align}
We instead invoke a result of \cite{leppanen2017} regarding convergence of conditional measures (Fact \ref{condec_sum} below), which implies an upper bound to \eqref{intro_approx6} depending on \(I\). Together with Jensen's inequality the bound then leads to \eqref{intro_approx4}.

\subsection{Notations} For comparing quantities, the following notations are used: Given two real-valued functions \(f\) and \(g\), we denote \(g(x) \lesssim_{\theta} f(x)\) if there exists a constant \(C > 0\) depending on \(\theta\) with \(g(x) \le Cf(x)\). Moreover, \(\lesssim \) means \(\lesssim_{\beta_*}\), i.e. that the constant depends only on the system \(T_{\beta_*}\).\\
\indent Given a measurable map \(f: [0,1] \to \bR\), we denote \(\Vert f \Vert_{\infty} = \sup_{x \in [0,1]} |f(x)| \) and \( \Vert f \Vert_{\text{Lip}} = \Vert f \Vert_{\infty} + \text{Lip}(f)\), where
\begin{align*}
\text{Lip}(f) = \sup_{x \neq y} \frac{|f(x) - f(y)|}{|x-y|}.
\end{align*}
If \(\mu\) is any Borel measure on \([0,1]\) such that \(f\) is \(\mu\)-integrable, we denote \(\mu(f) = \int f \, d \mu \). The Lebesgue measure on \([0,1]\) is denoted by \(m\).

\subsection*{Acknowledgements} I wish to thank Romain Aimino for several valuable discussions during the early stages of the paper's preparation. I would also like to thank my advisor Mikko Stenlund for suggesting the problem and for giving comments on preliminary versions of the manuscript. I gratefully acknowledge financial support from the Jane and Aatos Erkko Foundation, the Emil Aaltosen Säätiö, and Domast.

\section{Preliminaries}\label{prelim}

Throughout this section we consider a general time-dependent sequence \((T_{\alpha_n})_{n \ge 1} \) of Pomeau-Manneville maps, where \((\alpha_n)_{n\ge 1}\) is a sequence of numbers with \(0 \le \alpha_n \le \beta_*\) for some fixed \(\beta_* \in (0,1)\). We call such sequences of maps admissible. We  abbreviate \(T_{\alpha_n} = T_n\), \(n \ge 1\), and \(\cC_* = \cC_*(\beta_*)\), where \(\cC_*(\beta_*)\) is the cone of functions defined in the previous section (see below Definition \ref{def_int_qds}).

For each \(\alpha \in [0,1]\), the transfer operator associated to \(T_{\alpha}\) is denoted by \(\cL_{\alpha}\):
\begin{align*}
\cL_{\alpha}f(x) = \sum_{y \in T_{\alpha}^{-1}x} \frac{f(y)}{T_{\alpha}'(y)}, \hspace{0.5cm} f \in L^{1}([0,1],\bR).
\end{align*}

In the case of a single map \(T_{\alpha}\), we denote \(T_{\alpha}^{n+1} = T_{\alpha}^{n} \circ T_{\alpha}\) and \(\cL_{\alpha}^{n+1} = \cL_{\alpha}^{n} \cL_{\alpha}\), where \(T_{\alpha}^{0} = \text{id}_{[0,1]}\) and \(\cL_{\alpha}^{0} = \text{id}_{L^1}\). In the time-dependent setting, we introduce for all integers \(n \ge m\) the following notations:
\begin{align*}
T_n &= T_{\alpha_n} \hspace{3cm} \cL_n = \cL_{\alpha_n} \\
\widetilde{T}_{n,m} &= T_n \circ \cdots \circ T_m \hspace{1cm} \widetilde{\cL}_{n,m} = {\cL}_n  \cdots{\cL}_m \\
\widetilde{T}_n &= \widetilde{T}_{n,1} \hspace{2.95cm} \widetilde{\cL}_n = \widetilde{\cL}_{n,1}
\end{align*}
The map \(\widetilde{T}_{n,m}\) has \(2^{n-m +1} \) branches, and we denote the leftmost branch by \((\widetilde{T}_{n,m})_1\). The domain of the map \((\widetilde{T}_{n,m})_1\) is an interval whose left endpoint is the origin.

\subsection{On time-dependent intermittent maps.} Statistical properties of time-dependent intermittent systems have been studied before. The authors of \cite{bahsoun2014,bahsoun2016,bahsoun2017, bahsoun2017_2} have obtained various limit theorems in the setting of random intermittent maps, while in \cite{nicol2016} central limit theorems for sequential and random intermittent systems were shown. The existence of extreme value laws was proved recently in \cite{freitas2016}. For the present manuscript, \cite{aimino2015,leppanen2016,leppanen2017} are most important. We briefly review some results of these three papers related to correlation decay and perturbation of transfer operators.

\indent Let \((T_n)_{n\ge1}\) be an admissible sequence of maps. In \cite{aimino2015}, Aimino et al. proved polynomial memory loss for the sequential system.

\begin{fact}\label{aimino}
Let $f,g\in\cC_*$ with $\int f\,dx = \int g\,dx$. Then, for all integers $n\ge 0$,
\beqn
\|\widetilde\cL_n (f-g)\|_1 \lesssim (\|f\|_1+\|g\|_1) \rho(n),
\eeqn
where \(\rho(n) = n^{-\frac{1}{\beta_*}+1}(\log n)^{\frac{1}{\beta_*}}\) for \(n \ge 2\), and \(\rho(0) = \rho(1) = 1\).
\end{fact}

The proof of Fact \ref{aimino} was based on the earlier work ~\cite{liverani1999} where a similar result was obtained in the setting of a single map instead of a sequence of maps. \\
\indent In \cite{leppanen2017}, we observed that the method of \cite{aimino2015,liverani1999} can be extended to obtain a version of Fact \ref{aimino} for conditional densities. Given \(n \ge 1\), there is a partition \(\cP = \{ I_{\theta}^n \}_{\theta=1}^{2^n} \) of \((0,1)\) into open subintervals \(I_{\theta}^n\) such that \(\widetilde{T}_n \upharpoonright I_{\theta}^n\) maps \(I^n_{\theta}\) one-to-one and onto \((0,1)\) for all \(\theta \in \{ 1,\ldots, 2^n \}\). Given a probability measure \(\mu\) with density \(h \in \cC_*\), we define the conditional densities
\begin{align*}
h_{\theta} = \mu(I_{\theta}^n)^{-1} \mathbf{1}_{I_{\theta}^n} h, \hspace{0.5cm} \theta \in \{1,\ldots,2^n\},
\end{align*}
and denote \( \widetilde{h}_{\theta} = \cL_{n} \cdots \cL_{1}(h_{\theta}) = \widetilde{\cL}_{n}(h_{\theta})\).

\begin{fact}\label{condec_sum} Let \(h,g \in \cC_*\) be densities, and let \(m \ge 1\) be an integer. Then, 
\begin{align*}
& \sum_{\theta=1}^{2^n} \mu(I_{\theta}^n)\Vert \widetilde{\cL}_{n + m,n+1}(\widetilde{h}_{\theta} - \widetilde{\cL}_{n}g) \Vert_1 \lesssim \rho (m).
\end{align*}
 \end{fact}
Fact \ref{condec_sum} was proved in \cite{leppanen2017} and applied  together with Fact \ref{aimino} to show the following multicorrelation bound.
 
\begin{fact}\label{multicor} Let \(f_0,\ldots, f_{l}\) be Lipschitz continuous functions \([0,1] \to \bR\), and \(f_{l+1},\ldots,f_k \in L^{\infty}([0,1])\). Fix integers \(0 = n_0 \le n_1 \le \ldots \le n_k\), and denote
\begin{align*}
H &= f_0 \cdot f_1 \circ \widetilde{T}_{n_1} \cdots f_{l}\circ \widetilde{T}_{n_{l}} \\
G &=  f_{l+1} \circ \widetilde{T}_{n_{{l}+1}} \cdots f_{k}\circ \widetilde{T}_{n_{k}}.
\end{align*}
Then, for any probability measure \(\mu\) with density \(h \in \cC_*\),
\begin{align*}
&\left| \int HG \, d\mu - \int H \, d\mu \int G \, d\mu \right| \notag\\ 
&\lesssim  \prod_{i=0}^l \Vert f_i \Vert_{\textnormal{Lip}}\prod_{i=l+1}^k \Vert f_i \Vert_{\infty} \rho(n_{l+1}-n_{l}), 
\end{align*}
where \(\rho(n) = n^{-\frac{1}{\beta_*}+1}(\log n)^{\frac{1}{\beta_*}}\) for \(n \ge 2\), and \(\rho(0) = \rho(1) = 1\). \end{fact}

We end this review by recalling a couple of perturbation estimates from \cite{leppanen2016}.

\begin{fact}\label{SRB_continuous} For all \(h \in \cC_*\), 
\begin{align*}
\|(\cL_\alpha - \cL_\beta)h\|_1 \lesssim \|h\|_1(\beta-\alpha)^{\frac13(1- \beta_*)}|{\log(\beta-\alpha)}| 
\end{align*}
and
\begin{align*}
\|\hat h_\alpha-\hat h_\beta\|_1 \lesssim (\beta-\alpha)^{\frac13(1- \beta_*)^2} |{\log (\beta-\alpha)}|^\frac1{\beta_*}
\end{align*}
hold whenever $0\le\alpha<\beta\le \beta_*$.
\end{fact}

\section{The process \(\xi_n\)}

In the rest of this paper we consider the intermittent QDS described in Definition \ref{def_int_qds}. Throughout we assume that \(\gamma : [0,1] \to [0,1)\) is Hölder-continuous of order \(\eta \in (0,1]\), and that there exists \(\beta_* \in (0,1)\) such that \(\gamma([0,1]) \subset [0,\beta_*]\) and
\begin{align}\label{gamma_as_2}
\lim_{n\to\infty} n^{\eta}\sup_{t \in [0,1]} | \alpha_{n,\round{nt}} - \gamma_t | < \infty.
\end{align}
We fix a Lipschitz continuous function \(f: [0,1] \to \bR\) and an initial probability measure \(\mu\) whose density \(h \in \cC_*\). Recall that \(f_{n,k} = f \circ T_{n,k} \circ \cdots \circ T_{n,1}\) for any integers \(0 \le k \le n\). We denote
\begin{align*}
\xi_n(x,t) = \chi^{\mu}_n(x,t) = n^{\frac12}\int_0^t \bar{f}_{n,\round{ns}}(x) \, ds,   
\end{align*}
where
\begin{align*}
\bar{f}_{n,k} = f - \mu(f_{n,k}), \hspace{0.5cm}  0 \le k \le n.
\end{align*}
That is, we center according to the initial measure \(\mu\). \\
\indent Given integers \(0 \le k \le n\), we denote \(h_{n,k} = \cL_{n,k}\cdots \cL_{n,1} h\). In other words, \(h_{n,k}\) is the density of the pushforward measure \((T_{n,k}\circ \cdots \circ T_{n,1})_*\mu\). If \(r\) is close to \(s\), the density \(h_{n,\round{nr}}\) is pretty close to the SRB density \(\hat{h}_{\gamma_s}\), provided that the systems has been running for a while:

\begin{lem}\label{decay1} There exist \(p_0\in (0,1)\), \(p_1 \in (0,\tfrac12)\) and \(c_1 > 0\), such that
\begin{align*}
\Vert h_{n,\round{nr}} - \hat{h}_{\gamma_s} \Vert_1 \lesssim_{\gamma} n^{-p_0} +  |r-s|^{\eta\frac14 (1-\beta_*)^2}
\end{align*}
holds whenever \(s \in [0,1]\) and \( c_1n^{p_1-1} < r \le 1\).
\end{lem}
\begin{myproof} Let
\begin{align*}
p_1 = \tfrac{ \frac{\eta}{4} (1 - \beta_*)}{\frac{1}{\beta_*} + \frac{\eta}{4} (1 - \beta_*)} < \tfrac{1}{2}.
\end{align*}
The proof of Lemma 6.1 in \cite{leppanen2016} shows that for some \(p_0 \in (0,1)\) and \(c_1 > 0\),
\begin{align*}
\Vert h_{n,\round{nr}} - \hat{h}_{\gamma_r} \Vert_1 \lesssim_{\gamma} n^{-p_0}
\end{align*}
holds whenever \( c_1n^{p_1-1} < r \le 1\). On the other hand, Fact \ref{SRB_continuous} implies
\begin{align*}
\Vert \hat{h}_{\gamma_r} - \hat{h}_{\gamma_s} \Vert_1 \lesssim |\gamma_r-\gamma_s|^{\frac14 (1-\beta_*)^2},
\end{align*}
and the desired bound now follows by Hölder continuity of \(\gamma\).
\end{myproof}

\subsection{The variance \(\hat{\sigma}_t^2\)}\label{variance} Throughout the rest of this paper we assume \(\beta_* < \tfrac12\).
Note that by \eqref{gamma_as_2} there exists \(\beta_{**} < \tfrac12\), such that \(\alpha_{n,k} \le \beta_{**}\) whenever \(n\) is large enough and \(0 \le k \le n\). Since we are interested only in the limit \(n\to\infty\), we will assume without loss of generality that \(\alpha_{n,k} \le  \beta_* < \tfrac12\) for all \(n\) and \(0 \le k \le n\). \\
\indent Recall that \(\hat{f}_t = f -  \hat{\mu}_{\gamma_t}(f)\), and
\begin{align}\label{var}
\hat{\sigma}^2_t(f) = \lim_{m\to \infty} \hat{\mu}_{\gamma_t} \left[ \left( \frac{1}{\sqrt{m}} \sum_{k=0}^{m-1} \hat{f}_t \circ T_{\gamma_t}^k \right)^2 \right].
\end{align}

\begin{lem}\label{contin} For all \(t \in [0,1]\), $\hat{\sigma}^2_t(f)$ is finite and 
\begin{align}\label{var2}
\hat{\sigma}^2_t(f) = \hat{\mu}_{\gamma_t}[\hat{f}^2_t] + 2 \sum_{k=1}^{\infty} \hat{\mu}_{\gamma_t}[\hat{f}_t \hat{f}_t \circ T_{\gamma_t}^k].
\end{align}
Moreover, the map \(t \mapsto \hat{\sigma}^2_t(f)\) is Hölder continuous.
\end{lem}

\begin{myproof} A straightforward manipulation yields
\begin{align*}
\hat{\mu}_{\gamma_t} \left[ \left( \frac{1}{\sqrt{m}} \sum_{k=0}^{m-1} \hat{f}_t \circ T_{\gamma_t}^k \right)^2 \right] = \hat{\mu}_{\gamma_t}[\hat{f}^2_t] + 2  m^{-1} \sum_{k=1}^{m-1} (m-k) \hat{\mu}_{\gamma_t}[\hat{f}_t \hat{f}_t \circ T_{\gamma_t}^k],
\end{align*}
for all $m \ge 1$. By Fact \ref{multicor},
\begin{align*}
|\hat{\mu}_{\gamma_t}[\hat{f}_t \hat{f}_t \circ T_{\gamma_t}^k]| \lesssim \Vert f \Vert_{\text{Lip}}^2 \rho(k),
\end{align*}
which implies \(\sum_{k=1}^{\infty}  |\hat{\mu}_{\gamma_t}[\hat{f}_t \hat{f}_t \circ T_{\gamma_t}^k]| < \infty\) and  \(\lim_{m\to\infty} \tfrac1m \sum_{k=1}^{m-1}  k|\hat{\mu}_{\gamma_t}[\hat{f}_t \hat{f}_t \circ T_{\gamma_t}^k]| = 0\), since \(\beta_* < \tfrac12\). It follows that $\hat{\sigma}^2_t(f)$ is finite and the representation \eqref{var2} holds. \\
\indent To show the last claim, we define
\begin{align*}
V_{K,t} = \hat{\mu}_{\gamma_t}[\hat{f}^2_t] + 2 \sum_{k=1}^{K} \hat{\mu}_{\gamma_t}[\hat{f}_t \hat{f}_t \circ T_{\gamma_t}^k].
\end{align*}
Then,
\begin{align*}
&|\hat{\sigma}^2_t(f) - \hat{\sigma}^2_s(f)| 
\le  | V_{K,t} - V_{K,s} | + 2\sum_{k = K +1}^{\infty} \sup_r |\hat{\mu}_{\gamma_r}[\hat{f}_r \hat{f}_r \circ T_{\gamma_r}^k]| \\
&\lesssim_f  | V_{K,t} - V_{K,s} | + K^{1- \kappa},
\end{align*}
for some \(\kappa > 1\). To obtain a bound for \(| V_{K,t} - V_{K,s} |\), we first note that
\begin{align*}
V_{K,t} = [\hat{\mu}_{\gamma_t}(f^2) - \hat{\mu}_{\gamma_t}(f)^2] + \sum_{k=1}^K [\hat{\mu}_{\gamma_t}(f f \circ T_{\gamma_t}^k) - \hat{\mu}_{\gamma_t}(f)^2].
\end{align*}
Hence,
\begin{align*}
| V_{K,t} - V_{K,s} | 
&\le |\hat{\mu}_{\gamma_t}(f^2) - \hat{\mu}_{\gamma_s}(f^2) | + |\hat{\mu}_{\gamma_t}(f)^2 - \hat{\mu}_{\gamma_s}(f)^2| \\
&+ \sum_{k=1}^K |\hat{\mu}_{\gamma_t}(f f \circ T_{\gamma_t}^k) - \hat{\mu}_{\gamma_s}(f f \circ T_{\gamma_s}^k)| + \sum_{k=1}^K |\hat{\mu}_{\gamma_t}(f)^2 - \hat{\mu}_{\gamma_s}(f)^2|.
\end{align*}
By Fact \ref{SRB_continuous}, there is \(\kappa_1 \in (0,1)\) such that
\begin{align}\label{pert1}
\Vert (\cL_{\gamma_t} - \cL_{\gamma_s} )g \Vert_1 \lesssim \Vert g \Vert_1 |\gamma_t - \gamma_s |^{\kappa_1} \hspace{0.5cm} \forall g \in \cC_*,
\end{align}
and
\begin{align}\label{pert2}
|\hat{\mu}_{\gamma_t}(g) - \hat{\mu}_{\gamma_s}(g) | \lesssim \Vert g \Vert_{\infty} |\gamma_t - \gamma_s |^{\kappa_1} \hspace{0.5cm} \forall g \in L^{\infty}([0,1],\bR).
\end{align}
Since \(\gamma\) is Hölder continuous of order \(\eta\), it follows that \(\eqref{pert1}, \eqref{pert2} \lesssim_{g,\gamma} |t - s |^{\kappa_2} \), where \(\kappa_2 = \kappa_1 \eta \). Consequently,
\begin{align*}
|\hat{\mu}_{\gamma_t}(f^2) - \hat{\mu}_{\gamma_s}(f^2) | + |\hat{\mu}_{\gamma_t}(f)^2 - \hat{\mu}_{\gamma_s}(f)^2| \lesssim_{f,\gamma} |t - s |^{\kappa_2}.
\end{align*}
Moreover,
\begin{align*}
&|\hat{\mu}_{\gamma_t}(f f \circ T_{\gamma_t}^k) - \hat{\mu}_{\gamma_s}(f f \circ T_{\gamma_s}^k)| \\
&\le |\hat{\mu}_{\gamma_t}(f f \circ T_{\gamma_t}^k) - \hat{\mu}_{\gamma_s}(f f \circ T_{\gamma_t}^k)| + |\hat{\mu}_{\gamma_s}(f f \circ T_{\gamma_t}^k) - \hat{\mu}_{\gamma_s}(f f \circ T_{\gamma_s}^k) | \\
&\lesssim_{f,\gamma} |t - s |^{\kappa_2} + \Vert (\cL^k_{\gamma_t} - \cL^k_{\gamma_t})(f \hat{h}_{\gamma_s}) \Vert_1.
\end{align*}
Since
\begin{align*}
(\cL^k_{\gamma_t} - \cL^k_{\gamma_s})(f \hat{h}_{\gamma_s}) = \sum_{m =1}^k \cL_{\gamma_t}^{m-1} (\cL_{\gamma_t} - \cL_{\gamma_s})\cL_{\gamma_s}^{k-m} (f \hat{h}_{\gamma_s}),
\end{align*}
and \(\cL_{\gamma_t}\) is an \(L^1\)-contraction,
\begin{align*}
\Vert (\cL^k_{\gamma_t} - \cL^k_{\gamma_s})(f \hat{h}_{\gamma_s}) \Vert_1 \le 
\sum_{m =1}^k  \Vert (\cL_{\gamma_t} - \cL_{\gamma_s})\cL_{\gamma_s}^{k-m} (f \hat{h}_{\gamma_s}) \Vert_1.
\end{align*}
The function \(\cL_{\gamma_s}^{k-m} (f \hat{h}_{\gamma_s})\) might not belong to the cone \(\cC_*\), so we can not directly apply \eqref{pert1}. However, by Lemma 2.4 in \cite{leppanen2017} there exist \(g_1, \ldots, g_4 \in \cC_*\) such that
\begin{align*}
\cL_{s}^{k-m} (f \hat{h}_s) = g_1 - g_2 + g_3 - g_4,
\end{align*}
and for all \(1 \le i \le 4\),
\begin{align*}
\Vert g_i \Vert_1 \lesssim \Vert f \Vert_{\text{Lip}}.
\end{align*}
Combining this with \eqref{pert1} now yields the upper bound
\begin{align*}
\Vert (\cL^k_{\gamma_t} - \cL^k_{\gamma_s})(f \hat{h}_{\gamma_s}) \Vert_1 \lesssim_{f,\gamma} k|t - s |^{\kappa_2}.
\end{align*}
\indent With the foregoing bounds we conclude that
\begin{align*}
&| V_{K,t} - V_{K,s} | \lesssim_{f,\gamma} K^2 |t - s |^{\kappa_2}.
\end{align*}
Hence,
\begin{align*}
|\hat{\sigma}^2_t(f) - \hat{\sigma}^2_s(f)| \lesssim_{f,\gamma} K^2 |t - s |^{\kappa_2} + K^{1-\kappa} 
\end{align*}
holds for all \(K \ge 1\), \(s,t \in [0,1]\). We fix \(s,t \in [0,1]\) with \(s \neq t\), and set \(K = \round{|t -s|^{-\frac{\kappa_2}{\kappa +1}}} \) so that \(K^2 |t - s |^{\kappa_2} \approx K^{1-\kappa}  \). It follows that
\begin{align*}
|\hat{\sigma}^2_t(f) - \hat{\sigma}^2_s(f)| \lesssim_{f,\gamma} |t-s|^{\kappa_2 \frac{\kappa -1}{\kappa +1}}.
\end{align*}
This shows that the map \(t \mapsto \hat{\sigma}^2_t(f)\) is Hölder continuous, and thus completes the proof of Lemma \ref{contin}.
 \end{myproof}

\subsection{The second moment \(\mu[[\xi_n(t + \delta) - \xi_n(t)]^2]\)}\label{second_moment}  
The following simple estimate shows that the second moment 
\begin{align*}
\mu[[\xi_n(t + \delta) - \xi_n(t)]^2] = \int_0^1 [\xi_n(x,t + \delta) - \xi_n(x,t)]^2 \, d\mu(x)
\end{align*}
is uniformly bounded.

\begin{lem}\label{ubound} Whenever \(0 \le t \le t + \delta \le 1\),
\begin{align*}
\mu[[\xi_n(t + \delta) - \xi_n(t)]^2] \lesssim \Vert f \Vert^2_{\textnormal{Lip}}\delta.
\end{align*}
 \end{lem}

\begin{myproof} We start by writing
\begin{align*}
\mu[[\xi_n(t+\delta) - \xi_n(t)]^2] = \mu \left[n \left( \int^{t+\delta}_t \bar{f}_{n,\round{nu}} \, du \right)^2  \right] 
= n \int_{t}^{t+\delta} \int_t^{t+\delta} \mu(\bar{f}_{n,\round{nu}}\bar{f}_{n,\round{nv}}) \, du \, dv.
\end{align*}
Since \(\mu(\bar{f}_{n,k})=0\), Fact \ref{multicor} implies the bound
\begin{align*}
&|\mu(\bar{f}_{n,\round{nu}}\bar{f}_{n,\round{nv}})| \lesssim \Vert f \Vert_{\text{Lip}}^2 \rho (n(v-u)),
\end{align*}
whenever \(0 \le u \le v \le 1\). Hence, if \(\delta \le 2n^{-1}\),
\begin{align*}
\mu[[\xi_n(t+\delta) - \xi_n(t)]^2] 
\lesssim n\Vert f \Vert_{\text{Lip}}^2 \int_{t}^{t+\delta} \int_t^{t+\delta} 1 \, du \, dv \lesssim \Vert f \Vert_{\text{Lip}}^2  \delta.
\end{align*}
On the other hand, if \(\delta > 2n^{-1}\), there is \(\kappa = \kappa (\beta_*) > 1\) such that
\begin{align*}
&n \int_{t}^{t+\delta} \int_t^{t+\delta} \mu(\bar{f}_{n,\round{nu}}\bar{f}_{n,\round{nv}}) \, du \, dv \\
&\lesssim \Vert f \Vert^2_{\text{Lip}} \delta + \Vert f \Vert^2_{\text{Lip}} n \int_{t + 2/n}^{t+\delta} \int_t^{v - 2/n} (n(v-u))^{-\kappa} \, du \, dv \\
&\lesssim \Vert f \Vert^2_{\text{Lip}}\delta.
\end{align*}
In both cases we arrive at the bound of Lemma \ref{ubound}.
\end{myproof}

Before investigating further properties of the process \((\xi_n)\), we introduce for brevity the following notations, given any integers \(0 \le k \le l \le n\):
\begin{align*}
\widetilde{T}_{n,l,k} &= T_{n,l} \circ \cdots \circ T_{n,k} \hspace{1cm} \widetilde{\cL}_{n,l,k} = {\cL}_{n,l}  \cdots{\cL}_{n,k} \\
\widetilde{T}_{n,k} &= \widetilde{T}_{n,k,1} \hspace{2.95cm} \widetilde{\cL}_{n,k} = \widetilde{\cL}_{n,k,1}
\end{align*}

\begin{lem}\label{secmoment} Whenever \(0 \le t \le t + \delta \le 1\),
\begin{align*}
\mu[[\xi_n(t+\delta) - \xi_n(t)]^2] =  \int_{t}^{t+\delta} \hat{\sigma}^2_s(f) \, ds + \delta o(1) + o(n^{-\frac12}),
\end{align*}
as \(n \to \infty\), where the error terms are uniform in \(t\) and \(\delta\). Moreover,
\begin{align*}
\int_t^{t+\delta} \hat{\sigma}^2_s(f) \, ds = \delta\hat{\sigma}^2_t(f) + o(\delta),
\end{align*}
as \(\delta \to 0\), where the error term is uniform in \(t\).
  \end{lem}

\begin{myproof} The second claim follows immediately by uniform continuity of the function \(t \mapsto \hat{\sigma}^2_t(f)\) (Lemma \ref{contin}). We show the first claim. To this end, let \(\delta > 0\) and let \(n\) be sufficiently large so that \(n^{-\frac12} \le \delta\). Let \(p_1 \in (0,\tfrac12)\) and \(c_1\) be as in Lemma \ref{decay1}. We fix \(\kappa, \lambda \in (0,1)\) such that \(0 < \kappa < p_1 < \lambda < \frac12 \),
 and denote \(a_n = n^{-1+ \kappa}\), and \(b_n = n^{-1+\lambda}\). Then, we write
\begin{align*}
\mu[[\xi_n(t+\delta) - \xi_n(t)]^2] = n \int_{t}^{t+\delta} \int_t^{t+\delta} \mu(\bar{f}_{n,\round{nu}}\bar{f}_{n,\round{nv}}) \, du \, dv,
\end{align*}
and partition the domain of integration \([t,t+\delta]^2 = P_n \cup Q_n \cup R_n \ \), where
\begin{align*}
P_n = \{ (s,r) \in [t,t+\delta]^2 \: : \: t + b_n \le s \le t + \delta - b_n \text{ and } |r-s| \le a_n \},
\end{align*}
\begin{align*}
Q_n = \{ (s,r) \in [t,t+\delta]^2 \: : \: |r-s| \le a_n \text{ and either } s < t + b_n \text{ or } s > t + \delta - b_n \},
\end{align*}
and
\begin{align*}
R_n = \{ (s,r) \in [t,t+\delta]^2 \: : \: |r-s| > a_n \}.
\end{align*}
Since \(m(Q_n) = O(a_nb_n)\),
\begin{align*}
\left| n\iint_{Q_n} \mu(\bar{f}_{n,\round{ns}}\bar{f}_{n,\round{nr}}) \, dr \,ds \right| \lesssim \Vert f \Vert_{\infty}^2 na_nb_n \lesssim \Vert f \Vert_{\infty}^2 n^{-1 + \lambda + \kappa}.
\end{align*}
On the other hand, when \(n^{\kappa} = na_n \ge 2\), it follows by Fact \ref{multicor} that
\begin{align*}
&\left| n \iint_{R_n} \mu(\bar{f}_{n,\round{ns}}\bar{f}_{n,\round{nr}}) \, dr \,ds \right| = 
\left| 2n \int_t^{t+\delta-a_n} \int_{s+ a_n}^{t+\delta} \mu(\bar{f}_{n,\round{ns}}\bar{f}_{n,\round{nr}}) \, dr \,ds \right| \\
&\lesssim \Vert f \Vert_{\text{Lip}}^2 2n\int_t^{t+\delta-a_n} \int_{s+ a_n}^{t+\delta} \rho(nr-ns) \, dr \, ds  \\
&\lesssim \Vert f \Vert_{\text{Lip}}^2   2n\int_t^{t+\delta-a_n} \int_{s+ a_n}^{t+\delta} (nr-ns)^{-\theta} \, dr \, ds  \lesssim \Vert f \Vert_{\text{Lip}}^2\delta n^{- \kappa (\theta -1)},
\end{align*}
for some \(\theta = \theta(\beta_*) > 1\). Thus, if \(\kappa > 0\) is sufficiently small, all significant contribution comes from the diagonal \(P_n\):
\begin{align*}
\mu[[\xi_n(t+\delta) - \xi_n(t)]^2] = n \iint_{P_n} \mu(\bar{f}_{n,\round{ns}}\bar{f}_{n,\round{nr}}) \, dr \,ds + O(\delta n^{\kappa(1-\theta)} + n^{-1 + \lambda + \kappa}).
\end{align*}

Note that for all \((s,r) \in P_n\), we have the lower bound
\begin{align*} r &\ge s - a_n \ge t + b_n - a_n \ge b_n - a_n = (n^{\lambda-p_1} - n^{\kappa - p_1})n^{-1+p_1} \\
&\ge (n^{\lambda-p_1} - 1)n^{-1+p_1} > c_1n^{-1+p_1},
\end{align*} 
when \(n^{\lambda-p_1} > 1 + c_1\). For such \(n\), it follows by Lemma  \ref{decay1} that 
\begin{align}\label{estim1}
&\Vert h_{n,\round{nr}} - \hat{h}_{\gamma_s} \Vert_1 \lesssim_{\gamma} n^{-p_0} +  |r-s|^{\eta\frac14 (1-\beta_*)^2} 
\lesssim n^{-p_0} +  a_n^{\eta\frac14 (1-\beta_*)^2} \lesssim_{\gamma}  n^{-\theta_1},
\end{align}
for some \(\theta_1 = \theta_1(\beta_*,\eta) \in (0,1)\) (independent of \(\kappa\)), and all \((s,r) \in P_n\). Hence, 
\begin{align*}
&\sup_{r \in (s-a_n, s+ a_n)} | \mu(f_{n\round{nr}}) - \hat{\mu}_{\gamma_s}(f) |  \lesssim_{\gamma} \Vert f \Vert_{\infty} n^{-\theta_1},
\end{align*}
from which it follows that
\begin{align*}
& n \int_{s - a_n}^{s+a_n} \mu(\bar{f}_{n,\round{ns}}\bar{f}_{n,\round{nr}}) \, dr 
= n \int_{s - a_n}^{s+a_n} \mu(f_{n,\round{ns}} f_{n,\round{nr}}) - \hat{\mu}_{\gamma_s}(f)^2 \, dr  
+ O(n^{\kappa-\theta_1}).
\end{align*}
Next, we split the domain of integration \([s-a_n, s+a_n] = [s-a_n,s] \cup [s,s+a_n] \), and consider the right half. Letting \(c_n = \frac1n (1 - \{ns\}) \), we have
\begin{align*}
&n \int_s^{s+a_n} \mu(f_{n,\round{ns}} f_{n,\round{nr}}) \, dr = n \int_0^{a_n} \mu(f_{n,\round{ns}} f_{n,\round{n(s+r)}}) \, dr \\
&= c_n n \mu_{n,\round{ns}}(f^2) + n \int_{c_n}^{a_n} \mu_{n,\round{ns}}(f f\circ T_{n,\round{n(s+r)}} \circ \cdots \circ T_{n,\round{ns} +1} ) \, dr.
\end{align*}
Thus, by \eqref{estim1},
\begin{align*}
& n \int_s^{s+a_n} \mu(f_{n,\round{ns}} f_{n,\round{nr}}) \, dr \\
&= c_n n \hat{\mu}_{\gamma_s}(f^2) + n \int_{c_n}^{a_n} \hat{\mu}_{\gamma_s}(f f\circ T_{n,\round{n(s+r)}} \circ \cdots \circ T_{n,\round{ns} +1}) \, dr + O(n^{-\theta_1} + n^{\kappa-\theta_1}) \\
&=  n \int_0^{c_n} m(f\hat{h}_{\gamma_s}f) \, dr +  n \int_{c_n}^{a_n} m(f  \cL_{n,\round{n(s+r)}} \cdots \cL_{n,\round{ns} +1}(\hat{h}_{\gamma_s}f)) \, dr + O(n^{-\theta_1} + n^{\kappa-\theta_1}).
\end{align*}
We replace \(\widetilde{\cL}_{n,\round{n(s+r)},\round{ns}+1} = \cL_{n,\round{n(s+r)}} \cdots \cL_{n,\round{ns} +1}\) by \(\cL_{\gamma_s}^{\round{n(s+r)} - \round{ns}}\) in the second integral. To see that the inflicted error is small, first note that
\begin{align*}
&\Vert ( \widetilde{\cL}_{n,\round{n(s+r)},\round{ns}+1} -  \cL_{\gamma_s}^{\round{n(s+r)} - \round{ns}} ) \hat{h}_{\gamma_s} f \Vert_1 \\
&\le \sum_{k = \round{ns} + 1}^{\round{n(s+r)}} \Vert \widetilde{\cL}_{n,\round{n(s+r)},k+1} (\cL_{n,k} - \cL_{\gamma_s})\cL_{\gamma_s}^{k - \round{ns} - 1} \hat{h}_{\gamma_s} f \Vert_1 \\
&\lesssim nr \max_{\round{ns} +1 \le k \le \round{n(s+r)}} \Vert (\cL_{n,k} - \cL_{\gamma_s})\cL_{\gamma_s}^{k - \round{ns} - 1} \hat{h}_{\gamma_s} f \Vert_1.
\end{align*}
By Lemma 2.4 in \cite{leppanen2017}, for each \(k\) with \(\round{ns} +1 \le k \le \round{n(s+r)}\), there are functions \(g_1,\ldots,g_4 \in \cC_*\) such that
\begin{align*}
\cL_{\gamma_s}^{k - \round{ns} - 1}\hat{h}_{\gamma_s} f = g_1 - g_2 + g_3 - g_4, 
\end{align*}
and \( \Vert g_i \Vert_1 \lesssim \Vert f \Vert_{\text{Lip}}\). Fact \ref{SRB_continuous} applies to cone functions, and yields the bound
\begin{align*}
 \Vert (\cL_{n,k} - \cL_{\gamma_s})g_i\Vert_1 \lesssim \Vert f \Vert_{\text{Lip}} |\alpha_{n,k} - \gamma_s|^{\frac14 (1-\beta_*)}, \hspace{0.5cm} 1 \le i \le 4,
\end{align*}
so that for all \(r \le a_n\) and \(\round{ns} +1 \le k \le \round{n(s+r)}\),
\begin{align*}
&\Vert (\cL_{n,k} - \cL_{\gamma_s})\cL_{\gamma_s}^{k - \round{ns} - 1}\hat{h}_{\gamma_s} f \Vert_1 \\
&\lesssim \Vert f \Vert_{\text{Lip}} |\alpha_{n,k} - \gamma_s|^{\frac14 (1-\beta_*)} \\
& \lesssim_{\gamma} \Vert f \Vert_{\text{Lip}} \left(n^{-\eta} + \left|s- \frac{k}{n} \right|^{\eta}\right)^{\frac14 (1-\beta_*)} \\
&\lesssim_{\gamma} \Vert f \Vert_{\text{Lip}} (n^{-\eta} + r^{\eta} )^{\frac14 (1-\beta_*)}
\lesssim_{\gamma} \Vert f \Vert_{\text{Lip}} n^{-\theta_2},
\end{align*}
for some \(\theta_2 = \theta_2(\beta_*, \eta) \in (0,1)\) (independent of \(\kappa\)). We obtain the bound 
\begin{align*}
\Vert ( \widetilde{\cL}_{n,\round{n(s+r)},\round{ns}+1} -  \cL_{\gamma_s}^{\round{n(s+r)} - \round{ns}} ) \hat{h}_{\gamma_s} f \Vert_1  \lesssim nr \Vert f \Vert_{\text{Lip}} n^{-\theta_2} \lesssim \Vert f \Vert_{\text{Lip}} n^{\kappa - \theta_2}.
\end{align*}

By the foregoing estimates, we conclude that
\begin{align*}
&n \int_s^{s+a_n} \mu(f_{n,\round{ns}} f_{n,\round{nr}}) \, dr \\
&= n \int_0^{c_n} m(f\hat{h}_{\gamma_s}f) \, dr +  n \int_{c_n}^{a_n} m(f\cL_{\gamma_s}^{\round{n(s+r)} - \round{ns}} \hat{h}_{\gamma_s}f ) \, dr + O(na_nn^{\kappa-\theta_2} + n^{\kappa-\theta_1}) \\
& = n \int_0^{a_n} m(f\cL_{\gamma_s}^{\round{n(s+r)} - \round{ns}} \hat{h}_{\gamma_s}f) \, dr + O(n^{2\kappa-\theta_2} + n^{\kappa - \theta_1}).
\end{align*}
Similarly one shows that
\begin{align*}
n \int_{s-a_n}^{s} \mu(f_{n,\round{ns}} f_{n,\round{nr}}) \, dr = n \int_{-a_n}^{0} m(f\cL_{\gamma_s}^{\round{ns} - \round{n(s+r)}}(\hat{h}_{\gamma_s}f)) \, dr + O(n^{2\kappa-\theta_2} + n^{\kappa - \theta_1}).
\end{align*}
Hence,
\begin{align*}
&n \int_{s - a_n}^{s+a_n} \mu(\bar{f}_{n,\round{ns}}\bar{f}_{n,\round{nr}}) \, dr \\
&= n \int_{- a_n}^{a_n} m(f\cL_{\gamma_s}^{|\round{ns} - \round{n(s+r)}|} \hat{h}_{\gamma_s}\hat{f}_s ) \, dr + O(n^{2\kappa-\theta_2} + n^{\kappa - \theta_1}) \\
&= n \int_{- \infty}^{\infty} m(f\cL_{\gamma_s}^{|\round{ns} - \round{n(s+r)}|} \hat{h}_{\gamma_s}\hat{f}_s) \, dr + O(n^{\kappa(1 - \theta)} + n^{2\kappa-\theta_2} + n^{\kappa - \theta_1}) \\
&= \hat{\sigma}^2_s(f) + O(n^{\kappa(1 - \theta)} + n^{2\kappa-\theta_2} + n^{\kappa - \theta_1}),
\end{align*}
where Fact \ref{multicor} was used in the second equality. We have shown that
\begin{align*}
&\mu[[\xi_n(t+\delta) - \xi_n(t)]^2] \\
&= n \int_{t+b_n}^{t+\delta-b_n}\int_{s - a_n}^{s+a_n} \mu(\bar{f}_{n,\round{ns}}\bar{f}_{n,\round{nr}}) \, dr \,ds + \delta O(n^{\kappa(1-\theta)}) + O(n^{-1 + \lambda + \kappa}) \\
&= \int_{t+b_n}^{t+\delta-b_n} \hat{\sigma}^2_s(f)  \,ds + \delta O(n^{\kappa (1 - \theta)} + n^{2\kappa-\theta_2} + n^{\kappa - \theta_1}) + \delta O(n^{\kappa(1-\theta)}) + O(n^{-1 + \lambda + \kappa}) \\
&=\int_t^{t + \delta} \hat{\sigma}^2_s(f)  \,ds + \delta O(n^{\kappa (1 - \theta)} + n^{2\kappa-\theta_2} + n^{\kappa - \theta_1}) + \delta O(n^{\kappa(1-\theta)}) + \delta O(n^{-\frac12 + \lambda + \kappa}).
\end{align*}
The first claim of Lemma \ref{secmoment} follows from this by choosing sufficiently small \(\kappa < p_1\) and \(\lambda > p_1\).
\end{myproof}

\subsection{Decorrelation for the process $\xi_n$.}  Let \(C_c^{\infty}(\bR)\) denote the collection of all functions \(A \in C^{\infty}(\bR)\) with compact support. Together with Lemma \ref{secmoment}, the following decorrelation result forms the technical core for the proof of Theorem \ref{main}.

\begin{lem}\label{7.3} Let \(A \in C_c^{\infty}(\bR) \) and \(q \in \{1,2\}\). Then, for any probability measure \(\mu\) with density \(h \in \cC_*\),
\begin{align*}
\mu[A(\xi_n(s))[\xi_n(t) - \xi_n(s)]^q] - \mu[A(\xi_n(s)]\mu[[\xi_n(t) - \xi_n(s)]^q] = o(1),
\end{align*}
as \(n \to \infty\),  whenever \(0 \le s \le t \le 1\). Here the error is uniform in \(t\) and \(s\).
 \end{lem}

Before going into the proof, we record some preliminary bounds which will be instrumental later on too. If \(n \ge 1\) is an integer and \(t \in (0,1)\), then we can form the partition \(\cP = \{ I_{n,t,\theta} \}_{\theta=1}^{2^{\round{nt}}} \) of \((0,1)\) into open subintervals \(I_{n,t,\theta}\) such that the map \(\widetilde{T}_{n,\round{nt}} \upharpoonright I_{n,t,\theta}\) is one-to-one and onto \((0,1)\). Note that if \(I_{n,t,1}\) denotes the interval whose left endpoint is zero, then \(m(I_{n,t,1}) \ge m(I_{n,t,\theta})\) holds for all \(\theta \in \{1, \ldots, 2^{\round{nt}} \}\). Hence, whenever \(\round{nt}-\round{ns} \ge 1\), and \(\theta \in \{1, \ldots, 2^{\round{nt}} \}\),
\begin{align*}
&m(\widetilde{T}_{n,\round{ns}}(I_{n,t,\theta})) \le m( (\widetilde{T}_{n,\round{nt},\round{ns}+1})_1^{-1}(0,1) ) \le m((T_{\beta_*})_1^{-(\round{nt}-\round{ns})}(0,1)),
\end{align*}
where the last inequality follows from the fact that $\alpha \le \beta_*$ implies  $T_{\beta_*} \le T_{\alpha}$. By Lemma 3.2 in \cite{liverani1999},
\begin{align*}
m((T_{\beta_*})_1^{-(\round{nt}-\round{ns})}(0,1)) 
&\lesssim (nt-ns)^{-\frac{1}{\beta_*}},
\end{align*}
and so,
\begin{align*}
m(\widetilde{T}_{n,\round{ns}}(I_{n,t,\theta})) \lesssim (nt-ns)^{-\frac{1}{\beta_*}}.
\end{align*}
It follows that the function \(x \to \xi_n(x,s)\) is almost constant on each partition element \(I_{n,t,\theta}\), whenever \(\round{nt} \ge \round{ns} + 1\): For all \(x,y \in I_{n,t,\theta}\),
\begin{align}\label{eq:xi}
| \xi_n(x,s) - \xi_n(y,s)| \notag &\le n^{\frac12} \int_0^s |f_{n,\round{nu}}(x) - f_{n,\round{nu}}(y)| \, du \notag \\
 &\lesssim \Vert f \Vert_{\text{Lip}} n^{\frac12} \int_0^s (nt - nu)^{-\frac{1}{\beta_*}} \, du \notag\\
&\lesssim \Vert f \Vert_{\text{Lip}} \frac{1}{n^{\frac{1}{\beta_*}-\frac12}\left(\frac{1}{\beta_*}-1\right)}(t-s)^{1-\frac{1}{\beta_*}}.
\end{align}
\indent Next, suppose that \(B_1,\ldots,B_m\) are bounded and Lipschitz continuous functions on \(\bR\) and \(0 \le t_1 < \ldots < t_m < t_m + 2/n \le t \le 1\). For each \(\theta \in \{1, \ldots, 2^{\round{nt}} \}\), we denote by \(\mu_{n,t,\theta}\) the conditional measure \(\mu(I_{n,t,\theta})^{-1}\mu(1_{I_{n,t,\theta}}\cdot)\). Then, it follows from \eqref{eq:xi} that 
\begin{align}
&| B_1(\xi_n(x,t_1))\cdots B_m(\xi_n(x,t_m)) - \mu_{n,t,\theta}[B_1(\xi_n(t_1))\cdots B_m(\xi_n(t_m))]| \notag\\
&\lesssim \prod_{1 \le k \le m} \Vert B_k \Vert_{\text{Lip}} \Vert f \Vert_{\text{Lip}} \, n^{-\frac{1}{\beta_*} + \frac12}\sum_{k=1}^m(t-t_k)^{1-\frac{1}{\beta_*}}, \hspace{0.5cm} \forall x \in I_{n,t,\theta}.\label{eq:estim1}
\end{align}

With these preparations, we can now show Lemma \ref{7.3}:

\noindent Proof of Lemma \ref{7.3}:  We only prove the case \(q=2\), and leave the similar but easier case \(q=1\) to the reader. We need to show that whenever \(0 \le s \le t \le 1\), 
\begin{align}\label{to_show_1}
\mu[A(\xi_n(s))[\xi_n(t) - \xi_n(s)]^2] - \mu[A(\xi_n(s)]\mu[[\xi_n(t) - \xi_n(s)]^2] = o(1),
\end{align}
as \(n \to \infty\). \\
\indent Suppose that \(s < t\) and let \(n \ge 0\) be large enough so that \(n(t-s) \ge 2\). Note that \eqref{eq:estim1} implies the bound
\begin{align*}
&| \mu[1_{I_{n,s+2/n,\theta}}A(\xi_n(s))[\xi_n(t) - \xi_n(s)]^2] - \mu_{n,s + 2/n,\theta}[A(\xi_n(s))]\mu[1_{I_{n,s+2/n,\theta}}[\xi_n(t) - \xi_n(s)]^2]| \\
& \le \sup_{x \in I_{n,s+2/n,\theta}} | A(\xi_n(x,s)) - \mu_{n,s + 2/n,\theta}[A(\xi_n(s))]| \mu[1_{I_{n,s+2/n,\theta}}[\xi_n(t) - \xi_n(s)]^2] \\
& \lesssim \Vert A \Vert_{\text{Lip}} \, n^{-\frac{1}{\beta_*} + \frac12} \left(\tfrac{2}{n}\right)^{1-\frac{1}{\beta_*}}\mu[1_{I_{n,s+2/n,\theta}}[\xi_n(t) - \xi_n(s)]^2] \\
&\lesssim_A n^{-\frac12}\mu[1_{I_{n,s+2/n,\theta}}[\xi_n(t) - \xi_n(s)]^2], \hspace{0.5cm} \forall \theta \in \{1,\ldots, 2^{\round{ns}+2}\}.
\end{align*}
Recall that by Lemma \ref{ubound}, \( \mu[[\xi_n(t) - \xi_n(s)]^2]\) is bounded uniformly in \(s,t\) and \(n\). We partition \( [0,1] = \bigcup_{\theta =1}^{2^{\round{ns}+2}} I_{n,s+2/n,\theta}\) and utilize the above bound as follows:
\begin{align*}
&\mu[A(\xi_n(s))[\xi_n(t) - \xi_n(s)]^2] \\
&= \sum_{\theta=1}^{2^{\round{ns}+2}} \mu[1_{I_{n,s+2/n,\theta}}A(\xi_n(s))[\xi_n(t) - \xi_n(s)]^2] \\
&= \sum_{\theta=1}^{2^{\round{ns}+2}} \mu_{n,s + 2/n,\theta}[A(\xi_n(s))]\mu[1_{I_{n,s+2/n,\theta}}[\xi_n(t) - \xi_n(s)]^2] + 
O(n^{-\frac12}\mu[[\xi_n(t) - \xi_n(s)]^2]) \\
&= \sum_{\theta=1}^{2^{\round{ns}+2}} \mu[1_{I_{n,s+2/n,\theta}}A(\xi_n(s))]\mu_{n,s+2/n,\theta}[[\xi_n(t) - \xi_n(s)]^2] + 
O(n^{-\frac12}),
\end{align*}
where the error is uniform in \(t\) and \(s\). \\
\indent To obtain \eqref{to_show_1}, it remains to show that
\begin{align}\label{to_bound_1}
\sum_{\theta=1}^{2^{\round{ns}+2}} \mu(I_{n,s+2/n,\theta})|\mu_{n,s+2/n,\theta}[[\xi_n(t) - \xi_n(s)]^2] - \mu[[\xi_n(t) - \xi_n(s)]^2]| = o(1).
\end{align}
In other words, we need to replace the conditional measures \(\mu_{n,s+2/n,\theta}\) with \(\mu\) in the integral 
\begin{align*}
\mu_{n,s+2/n,\theta}[[\xi_n(t) - \xi_n(s)]^2] = n\iint_{[s,t]^2} \mu_{n,s+2/n,\theta}(\bar{f}_{n,\round{nu}}\bar{f}_{n,\round{nv}}) \, du \,dv, 
\end{align*}
and control the error. To achieve the latter goal, we first remove a small region \([s,s+ n^{-p}]^2\) from the domain of integration \([s,t]^2\), and then apply Fact \ref{condec_sum} together with the uniform boundedness of \(\mu[[\xi_n(t) - \xi_n(s)]^2]\). \\
\indent Let \(p \in (\tfrac12,1)\). By Cauchy-Schwarz inequality,
\begin{align*}
&n\left| \mu_{n,s+2/n,\theta}\left[\left[ \int_{s + n^{-p}}^t \bar{f}_{n,\round{nu}} \, du \right]^2\right] - \mu_{n,s+2/n,\theta}\left[\left[ \int_{s}^t \bar{f}_{n,\round{nu}} \, du \right]^2\right] \right| \\
&\lesssim_f n \left| \mu_{n,s+2/n,\theta}\left[ \left( \int_{s + n^{-p}}^t \bar{f}_{n,\round{nu}} \, du \right)\left( \int^{s + n^{-p}}_s \bar{f}_{n,\round{nv}} \, dv \right) \right] \right| + n^{1-2p} \\
&\lesssim_f n\left[ \mu_{n,s+2/n, \theta}\left( \int_{s + n^{-p}}^t \bar{f}_{n,\round{nu}} \, du \right)^2 \right]^{\frac12} \left[ \mu_{n,s+2/n,\theta}\left( \int^{s + n^{-p}}_s \bar{f}_{n,\round{nv}} \, dv \right)^2 \right]^{\frac12} \\
&+ n^{1-2p} \\
&\lesssim_f n^{1-p}\left[ \mu_{n,s+2/n, \theta}\left( \int_{s + n^{-p}}^t \bar{f}_{n,\round{nu}} \, du \right)^2 \right]^{\frac12} + n^{1-2p}.
\end{align*}
The increased lower limit \(s + n^{-p}\) in the remaining integral allows us to replace the conditional measures with a small error. First note that
\begin{align*}
&\left| \left[ \mu_{n,s+2/n, \theta}\left( \int_{s + n^{-p}}^t \bar{f}_{n,\round{nu}} \, du \right)^2 \right]^{\frac12} - \left[ \mu\left( \int_{s + n^{-p}}^t \bar{f}_{n,\round{nu}} \, du \right)^2 \right]^{\frac12} \right| \\
&\le \left[ \int_{s+n^{-p}}^t \int_{s+n^{-p}}^t  |\mu_{n,s+2/n, \theta}(\bar{f}_{n,\round{nu}}\bar{f}_{n,\round{nv}}) - \mu(\bar{f}_{n,\round{nu}}\bar{f}_{n,\round{nv}}) | \, du \, dv \right]^{\frac12} \\
&= \left[ 2 \int_{s+n^{-p}}^t \int_{s+n^{-p}}^v  |\mu_{n,s+2/n, \theta}(\bar{f}_{n,\round{nu}}\bar{f}_{n,\round{nv}}) - \mu(\bar{f}_{n,\round{nu}}\bar{f}_{n,\round{nv}}) | \, du \, dv \right]^{\frac12} \\
&\lesssim_f \left[ 2 \int_{s+n^{-p}}^t \int_{s+n^{-p}}^v  \Vert \widetilde{\cL}_{\round{nu}, \round{ns} + 2 +1} ( \widetilde{h}_{n,s+2/n,\theta} - \widetilde{\cL}_{\round{ns} +2}h) \Vert_1 \, du \, dv \right]^{\frac12}.
\end{align*}
Now, after an application of Jensen's inequality, Fact \ref{condec_sum} yields the bound
\begin{align}
&\sum_{\theta} \mu(I_{n,s+2/n,\theta})  \left| \left[ \mu_{n,s+2/n, \theta}\left( \int_{s + n^{-p}}^t \bar{f}_{n,\round{nu}} \, du \right)^2 \right]^{\frac12} - \left[ \mu\left( \int_{s + n^{-p}}^t \bar{f}_{n,\round{nu}} \, du \right)^2 \right]^{\frac12} \right| \notag\\
&\lesssim_f \sum_{\theta} \mu(I_{n,s+2/n,\theta}) \left[  \int_{s+n^{-p}}^t \int_{s+n^{-p}}^v  \Vert \widetilde{\cL}_{\round{nu}, \round{ns} + 2 +1} ( \widetilde{h}_{n,s+2/n,\theta} - \widetilde{\cL}_{\round{ns} +2}h) \Vert_1 \, du \, dv \right]^{\frac12} \notag\\
&\lesssim_f \left[  \int_{s+n^{-p}}^t \int_{s+n^{-p}}^v  \sum_{\theta} \mu(I_{n,s+2/n,\theta})  \Vert \widetilde{\cL}_{\round{nu}, \round{ns} + 2 +1} ( \widetilde{h}_{n,s+2/n,\theta} - \widetilde{\cL}_{\round{ns} +2}h) \Vert_1 \, du \, dv \right]^{\frac12} \notag\\
&\lesssim_f \left[  \int_{s+n^{-p}}^t \int_{s+n^{-p}}^v  \rho(\round{nu} - \round{ns} -2 )   \, du \, dv \right]^{\frac12}\notag \\
&\lesssim_f n^{-\frac{\kappa}{2}} \left[  \int_{s+n^{-p}}^t \int_{s+n^{-p}}^v  (u - s - \tfrac{2}{n})^{-\kappa} \, du \, dv \right]^{\frac12} \notag \\
&\lesssim_f n^{-\frac{\kappa}{2} + p \frac{\kappa - 1}{2} } (t-s)^{\frac12}. \label{eq:lower_limit}
\end{align}
Here \(\kappa > 1\), because \(\beta_* < \tfrac12\). Having replaced the conditional measures, we apply Lemma \ref{ubound}:
\begin{align*}
\left[ \mu\left( \int_{s + n^{-p}}^t \bar{f}_{n,\round{nr}} \, dr \right)^2 \right]^{\frac12} \lesssim_f n^{-\frac12} (t-s)^{\frac12}.
\end{align*}
We have shown that
\begin{align*}
&\sum_{\theta} \mu(I_{n,s+2/n,\theta}) n\left| \mu_{n,s+2/n,\theta}\left[\left[ \int_{s + n^{-p}}^t \bar{f}_{n,\round{nr}} \, dr \right]^2\right] - \mu_{n,s+2/n,\theta}\left[\left[ \int_{s}^t \bar{f}_{n,\round{nr}} \, dr \right]^2\right] \right| \\
&\lesssim_f n^{1-p} n^{-\frac{\kappa}{2} + p \frac{\kappa - 1}{2} } (t-s)^{\frac12} + n^{\frac12 -p } (t-s)^{\frac12} + n^{1-2p}.
\end{align*}
Similarly,
\begin{align*}
&n\left| \mu\left[\left[ \int_{s + n^{-p}}^t \bar{f}_{n,\round{nu}} \, du \right]^2\right] - \mu\left[\left[ \int_{s}^t \bar{f}_{n,\round{nu}} \, du \right]^2\right] \right| \\
&\lesssim_f n \left| \mu\left[ \left( \int_{s + n^{-p}}^t \bar{f}_{n,\round{nu}} \, du \right)\left( \int^{s + n^{-p}}_s \bar{f}_{n,\round{nv}} \, dv \right) \right] \right| + n^{1-2p} \\
&\lesssim_f n^{1-p}\left[ \mu\left( \int_{s + n^{-p}}^t \bar{f}_{n,\round{nr}} \, dr \right)^2 \right]^{\frac12}  + n^{1-2p} \\
&\lesssim_f n^{\frac12 -p }(t-s)^{\frac12} + n^{1-2p}.
\end{align*}

Combining the previous two bounds yields
\begin{align}
&\sum_{\theta=1}^{2^{\round{ns}+2}} \mu(I_{n,s+2/n,\theta})|\mu_{n,s+2/n,\theta}[[\xi_n(t) - \xi_n(s)]^2] - \mu[[\xi_n(t) - \xi_n(s)]^2]| \notag\\
&\lesssim_f \sum_{\theta} \mu(I_{n,s+2/n,\theta}) n\left| \mu_{n,s+2/n,\theta}\left[\left[ \int_{s + n^{-p}}^t \bar{f}_{n,\round{nu}} \, du \right]^2\right] - \mu\left[\left[ \int_{s + n^{-p}}^t \bar{f}_{n,\round{nu}} \, du \right]^2\right] \right| \notag\\
&+ n^{1-p} n^{-\frac{\kappa}{2} + p \frac{\kappa - 1}{2} } (t-s)^{\frac12}  + n^{\frac12 -p } (t-s)^{\frac12} + n^{1-2p} \notag\\
&\lesssim_f  n n^{ -\kappa + p(\kappa - 1)} (t-s) +  n^{1 - \frac{\kappa}{2} + p\frac{\kappa-3}{2} } (t-s)^{\frac12}  + n^{\frac12 -p } (t-s)^{\frac12} + n^{1-2p} \notag\\
&\lesssim_f   n^{ -(1-p)(\kappa - 1)} (t-s)  + n^{\frac12 -p } (t-s)^{\frac12} + n^{1-2p}. \label{cond_bound_0}
\end{align}
Since \(\kappa > 1\), the obtained upper bound vanishes as \(n \to \infty\). The proof of Lemma \ref{7.3} is complete. \(\qed\)

The preceding upper bound \eqref{cond_bound_0} will be useful in what follows. We remark that it holds for all \(0 \le s \le t \le 1\) and \(\tfrac12 < p < 1\). 

\section{Proof of Theorem \ref{main} }\label{main_pf} We now proceed to prove Theorem \ref{main}, starting with the case \(\nu = \mu\).  We recall that \(\bP^{\mu}_{n} = \bP^{\mu,\mu}_{n} \) denotes the  distribution of the random element \(x \mapsto \xi_n(x,\cdot)\), given an initial measure \(\mu\) with density \(h \in \cC_*\). Expectation with respect to \(\bP^{\mu}_{n}\) is denoted by \(\bE^{\mu}_n\).

\subsection{The fourth moment \(\mu[[\xi_n(t + \delta) - \xi_n(t)]^4] \). }

\begin{lem}\label{tight} If \(\beta_* < \tfrac13\),
\begin{align*}
\mu[[\xi_n(t + \delta) - \xi_n(t)]^4] \lesssim \Vert f \Vert^4_{\textnormal{Lip}}\delta^2
\end{align*}
holds whenever \(0 \le t \le t + \delta \le 1\). In particular, the sequence of measures \((\bP_n^{\mu})_{n\ge 1}\) is tight. If \( \beta_* < \tfrac12\), 
\begin{align*}
\mu[[\xi_n(t + \delta) - \xi_n(t)]^4] \lesssim \Vert f \Vert^4_{\textnormal{Lip}} \delta^2 (n\delta)^{2-\kappa} 
\end{align*}
holds for some \(1 < \kappa < 2\).
 \end{lem}

\begin{myproof} We assume, without loss of generality, that \(\delta > \tfrac2n\). By symmetry,
\begin{align*}
&\mu[[\xi_n(t + \delta) - \xi_n(t)]^4]
= \mu \left[ \left( n^{\frac12}\int_{t}^{t+\delta} \bar{f}_{n,\round{nu}} \, du \right)^4  \right] \\
&= n^2 \int_{t}^{t+\delta} \int_{t}^{t+\delta} \int_{t}^{t+\delta} \int_{t}^{t+\delta}  \mu(\bar{f}_{n,\round{ns}}\bar{f}_{n,\round{nr}}\bar{f}_{n,\round{nv}}\bar{f}_{n,\round{nu}})\, du \, dv \, dr \, ds \\
&= 4! n^2 \int_{t}^{t +\delta} \int_{t}^{s} \int_{t}^{r} \int_{t}^{v}  \mu(\bar{f}_{n,\round{ns}}\bar{f}_{n,\round{nr}}\bar{f}_{n,\round{nv}}\bar{f}_{n,\round{nu}})\, du \, dv \, dr \, ds.
\end{align*}
 For \(u \le v \le r \le s\), we use Fact \ref{multicor} to obtain
\begin{align*}
&|\mu(\bar{f}_{n,\round{ns}}\bar{f}_{n,\round{nr}}\bar{f}_{n,\round{nu}}\bar{f}_{n,\round{nv}})| \\
&\lesssim \Vert f \Vert^4_{\text{Lip}} \min \{\rho(n(s-r)), \rho(n(v-u)) \} \\
&\lesssim \Vert f \Vert^4_{\text{Lip}} \rho(n(s-r))^{\frac12}\rho(n(v-u))^{\frac12}.
\end{align*}
It follows that
\begin{align*}
&\mu[[\xi_n(t + \delta) - \xi_n(t)]^4] 
\lesssim \Vert f \Vert^4_{\text{Lip}} n^2 4! \left( \int_{t}^{t+\delta} \int_{t}^v \rho(n(v-u))^{\frac12} \, du \, dv \right)^2.
\end{align*}
Since \(\beta_* < \tfrac12\), there is \(\kappa > 1\) such that
\begin{align}
&\int_{t}^{t+\delta} \int_{t}^v \rho(n(v-u))^{\frac12} \, du \, dv \notag\\
&\lesssim_f \frac{\delta}{n} + \int_{t + 2/n}^{t + \delta} \int_{t}^v \rho(n(v-u))^{\frac12} \, du \, dv \notag\\
&\lesssim_f \frac{\delta}{n} + n^{-\frac12 \kappa }\int_{t + 2/n}^{t+\delta} \int_{t}^{v-2/n} (v-u)^{-\frac12 \kappa} \, du \, dv. \label{tight_integral}
\end{align}
If \(\beta_* < \tfrac13\), we can choose \(\kappa > 2\) so that \(\eqref{tight_integral} \lesssim \delta n^{-1}\). If \(\tfrac13 \le \beta_* < \tfrac12\),  we have \(1 < \kappa < 2\) and \(\eqref{tight_integral} \lesssim \delta^2 (\delta n)^{-\frac{\kappa}{2}}\). Both bounds of Lemma \ref{tight} have been verified.
\end{myproof}

\subsection{The process \((M_t)\).} Following \cite{dobbs2016}, we define for each \(t \in [0,1]\) the evaluation functional \(\pi_t : \, C([0,1]) \to \bR\) by
\begin{align*}
\pi_t(\omega) = \omega(t),
\end{align*} 
and the differential operator
\begin{align*}
\mathscr{L}_t = \frac12 \hat{\sigma}^2_t(f) \frac{d^2}{dx^2},
\end{align*}
where we recall that \(f: \, [0,1] \to \bR\) is a Lipschtiz continuous function. Given a function \(A \in C_c^{\infty}(\bR) \), we define the stochastic process \((M_t)_{t \in [0,1]}\) by the formula
\begin{align*}
M_t = A \circ \pi_t - A \circ \pi_0 - \int^t_0 \mathscr{L}_s(A) \circ \pi_s \, ds.
\end{align*}
Note that
\begin{align*}
M_t(\xi_n(x,\cdot)) = A(\xi_n(x,t)) - A(\xi_n(x,0)) - \int_0^t \frac12\hat{\sigma}^2_s(f)A''(\xi_n(x,s)) \, ds. 
\end{align*}

The following result is at the heart of Theorem \ref{main}, for it implies that the weak limit of \((\xi_n)\), given that it exists, must be a diffusion with  \(\mathscr{L}_t\) as its generator.

\begin{lem}\label{mart} Assume that \(\bP\) is the weak limit of a subsequence \((\bP_{n_k}^{\mu})_{k \ge 1}\). Then, for any \(A \in C_c^{\infty}(\bR)\), the process \((M_t)_{t \in [0,1]}\) is a martingale with respect to \(\bP\) and the filtration \((\mathfrak{F}_t)_{t \in [0,1]}\), where \(\mathfrak{F}_t\) is the sigma-algebra on \(C([0,1])\) generated by \(\{ \pi_s \: : \: 0 \le s \le t\} \). \end{lem}

\begin{myproof} We denote expectation with respect to \(\bP\) by \(\bE = \bE_{\bP}\). Since \(A\) and \(A''\) are bounded, \(\bE(|M_t|) < \infty\) holds for all \(t \in [0,1]\). We need to show that
\begin{align*}
\bE[M_t - M_r \,|\, \mathfrak{F}_r] = 0, \hspace{0.5cm} 0 \le r < t \le 1,
\end{align*}
and this is equivalent to the condition
\begin{align}\label{eq:marting}
\bE[B_1\circ \pi_{t_1} \cdots B_m \circ \pi_{t_m}(M_t - M_r)] = 0,
\end{align}
whenever \(m \ge 1\), \(B_1,\ldots,B_m : \, \bR \to \bR\) are bounded Lipschitz continuous functions, and \(0 < t_1 < \ldots < t_m \le r < t \le 1\). Let us fix such numbers and functions. \\
\indent Let \(q \in (0,\frac12)\), and denote \(K_n = \round{n^q(t-r)}\) and \(\delta_n = (t-r)/K_n\). Then, for any integer \(n \ge 1\), we can write
\begin{align*}
M_t - M_r = \sum_{k=0}^{K_n -1} (M_{r + (k+1)\delta_n} - M_{r + k\delta_n}).
\end{align*}
Since \(\bP\) is the weak limit of \((\bP_{n_k}^{\mu})_{k \ge 1}\), \eqref{eq:marting} follows once we establish
\begin{align*}
\lim_{n\to\infty}K_n \sup_{r \le u \le t- \delta_n} |\mu[B_1(\xi_n(t_1))\cdots B_m(\xi_n(t_m))(M_{u + \delta_n} - M_u)(\xi_n)]| = 0.
\end{align*}
We proceed as in the proof of Lemma \ref{7.3} and first partition the unit interval into appropriately small subintervals. This enables us to utilize the fact that \(x \mapsto \xi_n(x,s)\) is nearly constant on each partition element. \\
\indent Let \(n \ge 1\) be an integer so large that \(r + 2/n < t\). We fix \(u \in [r,t-\delta_n]\), and recall that \(\cP = \{ I_{n,u + 2/n,\theta} \}_{\theta=1}^{2^{\round{nu}+2}} \) is the partition of \((0,1)\) into open subintervals \(I_{n,u + 2/n,\theta}\) such that for each \(\theta \in \{1,\ldots , 2^{\round{nu}+2}\}\) the map \(\widetilde{T}_{n,\round{nu}+2} \upharpoonright I_{n,u + 2/n,\theta}\) is one-to-one and onto \((0,1)\). Let \(x^*_{n,u,\theta}\) denote the midpoint of \(I_{n,u + 2/n,\theta}\), and set \(c_{n,u,\theta} = \xi_n(x^*_{n,u,\theta}, u)\). We define the function \(\xi^*_{n,u} : \, [0,1] \times [0,1] \to \bR\) by letting
\begin{align*}
\xi^*_{n,u}(x,s) =  \xi_n(x,s) - \xi_n(x,u) + c_{n,u,\theta} \hspace{0.5cm} \forall x \in I_{n,u + 2/n,\theta}.
\end{align*}
By \eqref{eq:xi}, for all \(s \in [0,1]\),
\begin{align*}
\sup_{x \in [0,1]} |\xi_n(x,s) - \xi^*_{n,u}(x,s)| = \max_{\theta} \sup_{x \in I_{n,u + 2/n,\theta}} |\xi_n(x,u) - \xi_n(x^*_{n,u,\theta}, u) | \lesssim \Vert f \Vert_{\text{Lip}}n^{-\frac12}.
\end{align*}
Consequently, 
\begin{align*}
\sup_{x \in [0,1]} |(M_{u + \delta_n} - M_u)(\xi_n(x,\cdot)) - (M_{u + \delta_n} - M_u)(\xi^*_{n,u}(x,\cdot))| \lesssim_A \Vert f \Vert_{\text{Lip}} n^{-\frac12}.
\end{align*}
Hence,
\begin{align*}
&\mu[B_1(\xi_n(t_1))\cdots B_m(\xi_n(t_m))(M_{u + \delta_n} - M_u)(\xi_n)] \\
&= \mu[B_1(\xi_n(t_1))\cdots B_m(\xi_n(t_m))(M_{u + \delta_n} - M_u)(\xi_{n,u}^*)] + O(n^{-\frac12}) \\
&= \sum_{\theta=1}^{2^{\round{nu}+2}} \mu[\mathbf{1}_{I_{n,u+2/n,\theta}}B_1(\xi_n(t_1))\cdots B_m(\xi_n(t_m))(M_{u + \delta_n} - M_u)(\xi_{n,u}^*)] + O(n^{-\frac12}).
\end{align*}
Since \(0 < t_1 < \ldots < t_m \le r \le u \le u + 2/n\), it follows by \eqref{eq:estim1} that
\begin{align*}
& \sum_{\theta=1}^{2^{\round{nu}+2}}\mu[\mathbf{1}_{I_{n,u+2/n,\theta}}B_1(\xi_n(t_1))\cdots B_m(\xi_n(t_m))(M_{u + \delta_n} - M_u)(\xi_{n,u}^*)] \\
&= \sum_{\theta=1}^{2^{\round{nu}+2}} \mu_{n,u + 2/n, \theta}[B_1(\xi_n(t_1))\cdots B_m(\xi_n(t_m))] \mu[\mathbf{1}_{I_{n,u+2/n,\theta}}(M_{u + \delta_n} - M_u)(\xi_{n,u}^*)] \\
&+ O(n^{-\frac12}) \\ 
&= \sum_{\theta=1}^{2^{\round{nu}+2}} \mu[\mathbf{1}_{I_{n,u+2/n,\theta}}B_1(\xi_n(t_1))\cdots B_m(\xi_n(t_m))] \mu_{n,u + 2/n, \theta}[(M_{u + \delta_n} - M_u)(\xi_{n,u}^*)] \\
&+ O(n^{-\frac12}).
\end{align*}
The error terms above are uniform in \(u\), and we have \(K_n O(n^{-\frac12}) = o(1)\). Moreover, the functions \(B_1,\ldots, B_m\) are bounded. Thus, it remains to show that
\begin{align}\label{mart_toshow}
\lim_{n\to \infty}\sup_{u \in [r,t-\delta_n]}K_n \sum_{\theta=1}^{2^{\round{nu}+2}} \mu(I_{n,u+2/n,\theta})| \mu_{n,u + 2/n, \theta}[(M_{u + \delta_n} - M_u)(\xi_{n,u}^*)]| = 0.
\end{align}
Note that
\begin{align*}
&(M_{u + \delta_n} - M_u)(\xi_{n,u}^*(x,\cdot)) \\
&= A(\xi_{n,u}^*(x,u + \delta_n)) - A(\xi_{n,u}^*(x,u)) - \int_u^{u + \delta_n} \mathscr{L}_s A(\xi^*_{n,u}(x,s)) \, ds.
\end{align*}
When \(x \in I_{n,u + 2/n, \theta}\), we can Taylor expand \(A(\xi_{n,u}^*(x,u + \delta_n))\) at \(\xi^*_{n,u}(x,u) = c_{n,u,\theta} \). By Taylor's theorem, there exists \(\kappa_{n,u,\theta}(x) \in \bR\) such that
\begin{align*}
&(M_{u + \delta_n} - M_u)(\xi_{n,u}^*(x,\cdot)) \\
&= A'(c_{n,u,\theta})[\xi_n(x,u + \delta_n) - \xi_n(x,u)] \\
& +\left[ \frac12 A''(c_{n,u,\theta})[\xi_n(x,u + \delta_n) - \xi_n(x,u)]^2 - \int_u^{u + \delta_n} \mathscr{L}_sA(\xi^*_{n,u}(x,s)) \, ds  \right] \\
&+ \frac16 A'''(\kappa_{n,u,\theta})[\xi_n(x,u + \delta_n) - \xi_n(x,u)]^3. 
\end{align*}
We will consider each of these three terms separately, and establish bounds which imply \eqref{mart_toshow} when \(q \in (0,\tfrac12)\) is chosen appropriately. We remark that everything we have done so far holds for arbitrary \(q \in (0,\tfrac12)\). \\
\indent Recall that \(h_{n,u+2/n,\theta}\) denotes the density of the conditional measure \(\mu_{n,u+2/n,\theta}\), and that 
\begin{align*}
\widetilde{h}_{n,u+2/n,\theta} = \cL_{n,\round{nu}+2} \cdots \cL_{n,1}h_{n,u+2/n,\theta}.
\end{align*}

\noindent\(\mathbf{1}^{\circ}\) The first term: We have
\begin{align*}
\mu_{n,u + 2/n,\theta}[A'(c_{n,u,\theta})[\xi_n(u + \delta_n) - \xi_n(u)]] = A'(c_{n,u,\theta})\mu_{n,u + 2/n,\theta}[\xi_n(u + \delta_n) - \xi_n(u)],
\end{align*}
where \(A'\) is bounded, and
\begin{align*}
&|\mu_{n,u + 2/n,\theta}[\xi_n(u + \delta_n) - \xi_n(u)]|\\
&\le n^{\frac12} \int^{u+\delta_n}_u |\mu_{n,u + 2/n,\theta}(f_{n,\round{ns}}) - \mu(f_{n,\round{ns}})| \, ds \\
&\lesssim_f n^{\frac12} \int^{u+\delta_n}_{u + 4/n} |\mu_{n,u + 2/n,\theta}(f_{n,\round{ns}}) - \mu(f_{n,\round{ns}})| \, ds + n^{-\frac12} \\
&\lesssim_f n^{\frac12} \int^{u+\delta_n}_{u + 4/n} \Vert \widetilde{\cL}_{n,\round{ns}, \round{nu} + 2 +1}(\widetilde{h}_{n,u+2/n,\theta} - \widetilde{\cL}_{n,\round{nu} + 2 }h) \Vert_1 \, ds + n^{-\frac12}.
\end{align*}
Fact \ref{condec_sum} can now be applied to obtain
\begin{align*}
&K_n \sum_{\theta=1}^{2^{\round{nu}+2}} \mu(I_{n,u+2/n,\theta})|\mu_{n,u + 2/n,\theta}[A'(c_{n,u,\theta})[\xi_n(u + \delta_n) - \xi_n(u)]]| \\
&\lesssim_{f,A} K_n \sum_{\theta=1}^{2^{\round{nu}+2}} \mu(I_{n,u+2/n,\theta})n^{\frac12} \int^{u+\delta_n}_{u + 4/n} \Vert \widetilde{\cL}_{n,\round{ns}, \round{nu} + 2 +1}(\widetilde{h}_{n,u+2/n,\theta} - \widetilde{\cL}_{n,\round{nu} + 2 }h) \Vert_1 \, ds \\
&+ K_n n^{-\frac12 } \\
&\lesssim_{A,f} K_n n^{\frac12} \int_{u + 4/n}^{u+\delta_n} \rho(\round{ns} - \round{nu}-2) \, ds + K_nn^{-\frac12} \\
&\lesssim_{A,f} \frac{1}{\kappa - 1} n^{q- \frac12} + n^{q- \frac12},
\end{align*}
where \(\kappa > 1\). Since \(q < \tfrac12\), \(n^{q- \frac12} = o(1)\).\\

\noindent\(\mathbf{2}^{\circ}\) The third term: Since \(A'''\) is bounded,
\begin{align*}
|\mu_{n,u + 2/n, \theta}[A'''(\kappa_{n,u,\theta})[\xi_n(u + \delta_n) - \xi_n(u)]^3]| \le \Vert A''' \Vert_{\infty} \mu_{n,u + 2/n, \theta}[|\xi_n(u + \delta_n) - \xi_n(u)|^3].
\end{align*}
We will bound the remaining third moment using essentially the same strategy as in the case of the second moment, which was treated in the proof of Lemma \ref{7.3}. This is feasible, since we control \(\mu[[\xi_n(u + \delta_n) - \xi_n(u)]^{\lambda}]\) for \(\lambda \le 4\).  \\
\indent Given \(p \in (\tfrac12, 1)\), Cauchy-Schwarz inequality yields
\begin{align*}
& \mu_{n,u + 2/n,\theta}\left[ \left| n^{\frac32} \left( \int_u^{u + \delta_n} \bar{f}_{n,\round{ns}} \, ds \right)^3  -  n^{\frac32} \left( \int_{u + n^{-p}}^{u + \delta_n}\bar{f}_{n,\round{ns}}\, ds \right)^3 \right| \right] \\
&\lesssim_f  n^{\frac32} \mu_{n,u + 2/n,\theta}\left[  \left( \int_u^{u + n^{-p}} \bar{f}_{n,\round{ns}}  \, ds \right)^2\left| \int_{u + n^{-p}}^{u + \delta_n} \bar{f}_{n,\round{ns}} \, ds \right| \right] \\
&+ n^{\frac32}\mu_{n,u + 2/n,\theta}\left[ \left| \int_u^{u + n^{-p}} \bar{f}_{n,\round{ns}} \, ds \right| \left( \int_{u + n^{-p}}^{u + \delta_n} \bar{f}_{n,\round{ns}}  \, ds \right)^2 \right]  \\
&+n^{\frac32 - 3p} \\
&\lesssim_f  n^{\frac32} \left[\mu_{n,u + 2/n,\theta}\left[ \left( \int_u^{u + n^{-p}} \bar{f}_{n,\round{ns}} \, ds \right)^4 \right] \right]^{\frac12} \left[ \mu_{n,u + 2/n,\theta}\left[\left( \int_{u + n^{-p}}^{u + \delta_n} \bar{f}_{n,\round{ns}} \, ds \right)^2 \right] \right]^{\frac12} \\
&+n^{\frac32} \left[\mu_{n,u + 2/n,\theta}\left[ \left( \int_u^{u + n^{-p}} \bar{f}_{n,\round{ns}} \, ds \right)^2 \right] \right]^{\frac12} \left[ \mu_{n,u + 2/n,\theta}\left[\left( \int_{u + n^{-p}}^{u + \delta_n} \bar{f}_{n,\round{ns}} \, ds \right)^4 \right] \right]^{\frac12} \\
&+n^{\frac32 - 3p} \\
&\lesssim_f n^{\frac32 -2p}   \left[ \mu_{n,u + 2/n,\theta}\left[\left( \int_{u + n^{-p}}^{u + \delta_n} \bar{f}_{n,\round{ns}} \, ds \right)^2 \right] \right]^{\frac12} \\
&+n^{\frac32-p}\left[ \mu_{n,u + 2/n,\theta}\left[\left( \int_{u + n^{-p}}^{u + \delta_n} \bar{f}_{n,\round{ns}} \, ds \right)^4 \right] \right]^{\frac12} +n^{\frac32 - 3p}.
\end{align*}
We replace both conditional measures \(\mu_{n,u + 2/n,\theta}\) with \(\mu\) in the remaining integrals. To control the error, we argue as in \eqref{eq:lower_limit} and observe that for \(\lambda \in \{2,4\}\) there is some \(\kappa > 1\) such that
\begin{align*}
&\sum_{\theta} \mu(I_{n,u+2/n,\theta}) \left| \left[ \mu_{n,u+2/n, \theta}\left( \int_{u + n^{-p}}^{u+\delta_n} \bar{f}_{n,\round{ns}} \, ds \right)^{\lambda} \right]^{\frac12} - \left[ \mu\left( \int_{u + n^{-p}}^{u+\delta_n} \bar{f}_{n,\round{ns}} \, ds \right)^{\lambda} \right]^{\frac12} \right| \\
&\lesssim_f n^{-\frac{\kappa}{2} + p \frac{\kappa - 1}{2} } \delta_n^{\frac{\lambda -1}{2}}.
\end{align*}
Then, Lemmas \ref{ubound} and \ref{tight} can be applied: 
\beqn
\left[ \mu\left[\left( \int_{u + n^{-p}}^{u + \delta_n} \bar{f}_{n,\round{ns}} \, ds \right)^{2} \right] \right]^{\frac12} \lesssim_f \delta_n^{\frac12} n^{-\frac12},
\eeqn
and, for some \(1 < \kappa_0 < 2\),
\begin{align*}
\left[ \mu\left[\left( \int_{u + n^{-p}}^{u + \delta_n} \bar{f}_{n,\round{ns}} \, ds \right)^{4} \right] \right]^{\frac12} \lesssim_f n^{-1} \delta_n (n\delta_n)^{\frac{2-\kappa_0}{2}}.
\end{align*}
These steps lead to the bound
\begin{align*}
& \sum_{\theta} \mu(I_{n,u+2/n,\theta}) \mu_{n,u + 2/n,\theta}\left[ \left| n^{\frac32} \left( \int_u^{u + \delta_n} \bar{f}_{n,\round{ns}} \, ds \right)^3  - n^{\frac32} \left( \int_{u + n^{-p}}^{u + \delta_n} \bar{f}_{n,\round{ns}} \, ds \right)^3 \right| \right] \\
&\lesssim_f n^{\frac32-2p} n^{-\frac{\kappa}{2} + p \frac{\kappa - 1}{2} } \delta_n^{\frac{2 -1}{2}} + n^{\frac32 - 2p} \delta_n^{\frac12} n^{-\frac12} \\
&+ n^{\frac32-p} n^{-\frac{\kappa}{2} + p \frac{\kappa - 1}{2} } \delta_n^{\frac{4 -1}{2}} + n^{\frac32 - p} n^{-1} \delta_n (n\delta_n)^{\frac{2-\kappa_0}{2}} + n^{\frac32-3p}  \\
&\lesssim_f 
 n^{1-2p}\delta_n^{\frac12} 
+ n^{\frac{(1-p)(3-\kappa)}{2}}\delta_n^{\frac32}  
+ n^{\frac12 - p} n^{ \frac{2-\kappa_0}{2}}  (\delta_n)^{1 + \frac{2-\kappa_0}{2}}  
+ n^{\frac32 - 3p} .
\end{align*}
By a similar (but easier) argument,
\begin{align*}
& \mu\left[ \left| n^{\frac32} \left( \int_u^{u + \delta_n} \bar{f}_{n,\round{ns}} \, ds \right)^3 -  n^{\frac32} \left( \int_{u + n^{-p}}^{u + \delta_n} \bar{f}_{n,\round{ns}} \, ds \right)^3  \right| \right] \\
&\lesssim_f   n^{1-2p}\delta_n^{\frac12} +  n^{\frac12 - p} n^{ \frac{2-\kappa_0}{2}}  (\delta_n)^{1 + \frac{2-\kappa_0}{2}}    + n^{\frac32 - 3p}.
\end{align*}
Hence,
\begin{align*}
&K_n \sum_{\theta} \mu(I_{n,u+2/n,\theta})\mu_{n,u + 2/n, \theta}[|[\xi_n(u + \delta_n) - \xi_n(u)]^3|]  \\
&\lesssim_f \delta_n^{-1} n^{\frac32} \sum_{\theta} \mu(I_{n,u+2/n,\theta}) \left| \mu_{n,u + 2/n,\theta}\left[  \left| \left( \int_{u + n^{-p}}^{u + \delta_n} \bar{f}_{n,\round{ns}} \, ds \right)^3 \right|\right] - \mu\left[ \left|\left( \int_{u + n^{-p}}^{u + \delta_n} \bar{f}_{n,\round{ns}} \, ds \right)^3 \right| \right] \right| \\
&+ n^{1-2p}\delta_n^{-\frac12} 
+ n^{\frac{(1-p)(3-\kappa)}{2}}\delta_n^{\frac12}  
+ n^{\frac12 - p} n^{ \frac{2-\kappa_0}{2}}  \delta_n^{\frac{2-\kappa_0}{2}}  
+ n^{\frac32 - 3p+q} \\
&\lesssim_f  \delta_n n^{ \frac12 -  (1-p)(\kappa-1)} + n^{1-2p}\delta_n^{-\frac12} 
+ n^{\frac{(1-p)(3-\kappa)}{2}}\delta_n^{\frac12}  
+ n^{\frac12 - p} (n\delta_n)^{\frac{2-\kappa_0}{2}}  
+ n^{\frac32 - 3p+q}, 
\end{align*}
where a bound similar to \eqref{eq:lower_limit} yields the second inequality. \\
\indent Recall that \(\delta_n \sim n^{-q}\), where \(0 < q < \tfrac12\). To guarantee that
\begin{align*}
\delta_n n^{ \frac12 -  (1-p)(\kappa-1)} + n^{1-2p}\delta_n^{-\frac12} 
+ n^{\frac{(1-p)(3-\kappa)}{2}}\delta_n^{\frac12}  
+ n^{\frac12 - p} (n\delta_n)^{\frac{2-\kappa_0}{2}}  
+ n^{\frac32 - 3p+q} = o(1),
\end{align*}
it suffices to choose
\begin{align*}
1 - \tfrac{q}{2} < p < 1 - \tfrac{\frac12 - q}{\kappa -1}.
\end{align*}
Such \(p\) exist, when \(q > \frac{1}{1+\kappa}\).\\

\noindent\(\mathbf{3}^{\circ} \) The second term. Recall that the remaining term we need to control is
\begin{align*}
 \frac12 A''(c_{n,u,\theta})[\xi_n(x,u + \delta_n) - \xi_n(x,u)]^2 - \int_u^{u + \delta_n} \mathscr{L}_sA(\xi^*_{n,u}(x,s)) \, ds.
\end{align*}
We have
\begin{align*}
&\mu_{n,u+2/n, \theta} \left[\frac12 A''(c_{n,u,\theta})[\xi_n(u + \delta_n) - \xi_n(u)]^2 \right]  \\
&=  \frac12 A''(c_{n,u,\theta}) \mu_{n,u+2/n, \theta} \left[ [\xi_n(u + \delta_n) - \xi_n(u)]^2 \right],
\end{align*}
and,
\begin{align*}
&\mu_{n, u + 2/n, \theta}\left[ \int_u^{u + \delta_n} \mathscr{L}_sA(\xi^*_{n,u}(s)) \, ds \right] \\
&= \int_{u}^{u+ \delta_n} \mu_{n, u + 2/n, \theta}\left[\frac12\hat{\sigma}^2_s(f)A''(\xi^*_{n,u}(s)) \right] \, ds \\
&= \frac12 \int_{u}^{u+ \delta_n} \hat{\sigma}^2_s(f)\mu_{n, u + 2/n, \theta}\left[A''(\xi^*_{n,u}(s)) \right] \, ds \\
&= \frac12 \int_{u}^{u+ \delta_n} \hat{\sigma}^2_s(f)\mu_{n, u + 2/n, \theta}\left[A''(\xi_n(s) - \xi_n(u) + c_{n,u,\theta}) \right] \, ds .
\end{align*}
By Taylor's theorem, there exists \(\widetilde{\kappa}_{n,u,\theta,s}(x) \in \bR\) such that
\begin{align*}
&\mu_{n, u + 2/n, \theta}\left[A''(\xi_n(s) - \xi_n(u) + c_{n,u,\theta}) \right] \\
&= A''(c_{n,u,\theta}) + \mu_{n, u + 2/n, \theta}\left[A'''(\widetilde{\kappa}_{n,u,\theta,s})(\xi_n(s) -\xi_n(u))\right].
\end{align*}
Let \(u \le s \le u + \delta_n\). Jensen's inequality combined with \eqref{cond_bound_0} implies the bound
\begin{align*}
&\sum_{\theta} \mu(I_{n,u+2/n,\theta})|\mu_{n,u+2/n,\theta} [[\xi_n(s) -\xi_n(u)]^2]^{\frac12} - \mu [[\xi_n(s) -\xi_n(u)]^2]^{\frac12}| \\
&\lesssim_f \left[ \sum_{\theta} \mu(I_{n,u+2/n,\theta})|\mu_{n,u+2/n,\theta} [[\xi_n(s) -\xi_n(u)]^2] - \mu [[\xi_n(s) -\xi_n(u)]^2]| \right]^{\frac12} \\
&\lesssim_f [ n^{ -(1-p)(\kappa - 1)} (s-u)  + n^{\frac12 -p } (s-u)^{\frac12} + n^{1-2p}]^{\frac12} \\
&\lesssim_f n^{ - \frac12 (1-p)(\kappa - 1)} (s-u)^{\frac12}  + n^{\frac12 (\frac12 -p) } (s-u)^{\frac14} + n^{\frac12 -p},
\end{align*}
where \(p \in (\tfrac12, 1)\) and \(\kappa > 1\). After an application of Cauchy-Schwarz inequality, we can apply the foregoing bound to obtain
\begin{align*}
&K_n \sum_{\theta} \mu(I_{n,u+2/n,\theta}) \frac12 \int_{u}^{u+ \delta_n} \hat{\sigma}^2_s(f) |\mu_{n, u + 2/n, \theta}\left[A'''(\widetilde{\kappa}_{n,u,\theta,s})(\xi_n(s) -\xi_n(u)) \right]| \, ds \\
&\lesssim_{f,A} K_n  \int_{u}^{u+ \delta_n} \sum_{\theta} \mu(I_{n,u+2/n,\theta}) |\mu_{n,u+2/n,\theta} [[\xi_n(s) -\xi_n(u)]^2]^{\frac12} - \mu [[\xi_n(s) -\xi_n(u)]^2]^{\frac12}| \, ds \\
&+K_n  \int_{u}^{u+ \delta_n} \mu [[\xi_n(s) -\xi_n(u)]^2]^{\frac12} \, ds \\
&\lesssim_{f,A} K_n \int_{u}^{u+\delta_n} n^{ - \frac12 (1-p)(\kappa - 1)} (s-u)^{\frac12}  + n^{\frac12 (\frac12 -p) } (s-u)^{\frac14} \, ds +  n^{\frac12 - p} + K_n \int_u^{u+\delta_n} (s-u)^{\frac12} \, ds \\
&\lesssim_{f,A} \delta_n^{\frac12} n^{ - \frac12 (1-p)(\kappa - 1)} + \delta_n^{\frac14} n^{\frac12 (\frac12 -p) } + n^{\frac12 - p} + \delta_n^{\frac12},
\end{align*}
where Lemma \ref{ubound} was used in the second inequality. The bound yields
\begin{align*}
&K_n \sum_{\theta} \mu(I_{n,u+2/n,\theta})\left|\mu_{n,u + 2/n, \theta}\left[ \frac12 A''(c_{n,u,\theta})[\xi_n(u + \delta_n) - \xi_n(u)]^2 - \int_u^{u + \delta_n} \mathscr{L}_sA(\xi^*_{n,u}(s)) \, ds  \right]\right| \\
&\lesssim_{f,A} K_n \sum_{\theta} \mu(I_{n,u+2/n,\theta})\left|\mu_{n,u + 2/n, \theta}\left[ [\xi_n(u + \delta_n) - \xi_n(u)]^2 - \int_{u}^{u+ \delta_n} \hat{\sigma}^2_s(f)  \, ds \right]\right|  \\
&+ K_n \sum_{\theta} \mu(I_{n,u+2/n,\theta}) \frac12 \int_{u}^{u+ \delta_n} \hat{\sigma}^2_s(f) |\mu_{n, u + 2/n, \theta}\left[A'''(\widetilde{\kappa}_{n,u,\theta,s})(\xi_n(s) -\xi_n(u)) \right]| \, ds  \\
&\lesssim_{f,A} K_n \sum_{\theta} \mu(I_{n,u+2/n,\theta})\left|\mu_{n,u + 2/n, \theta}\left[ [\xi_n(u + \delta_n) - \xi_n(u)]^2 - \int_{u}^{u+ \delta_n} \hat{\sigma}^2_s(f)  \, ds \right]\right| \\
&+\delta_n^{\frac14} n^{\frac12 (\frac12 -p) } + n^{\frac12 - p} + \delta_n^{\frac12},
\end{align*}
where \(\delta_n^{\frac14} n^{\frac12 (\frac12 -p) } + n^{\frac12 - p} + \delta_n^{\frac12} = o(1)\). \\
\indent It remains to  bound
\begin{align}\label{mart:final}
K_n \sum_{\theta} \mu(I_{n,u+2/n,\theta})\left|\mu_{n,u + 2/n, \theta}\left[ [\xi_n(u + \delta_n) - \xi_n(u)]^2 - \int_{u}^{u+ \delta_n} \hat{\sigma}^2_s(f)  \, ds \right]\right|.
\end{align}
To this end, we implement Lemma \ref{secmoment} together with \eqref{cond_bound_0} to obtain for any \(p \in (\tfrac12,1)\),
\begin{align*}
\eqref{mart:final} \le &K_n \sum_{\theta} \mu(I_{n,u+2/n,\theta}) \left|   \mu_{n,u+2/n, \theta}[ [\xi_n(u + \delta_n) - \xi_n(u)]^2]  -   \mu[[\xi_n(u + \delta_n) - \xi_n(u)]^2]\right| \\
&+ K_n \left|\mu[\xi_n(u + \delta_n) - \xi_n(u)]^2  -   \int_{u }^{u+ \delta_n} \hat{\sigma}^2_s(f) \, ds  \right| \\
&\lesssim_f K_n (n^{ -(1-p)(\kappa - 1)} \delta_n  + n^{\frac12 -p } \delta_n^{\frac12} + n^{1-2p}) + K_n \delta_n o(1) + n^{q - \frac12}\\
&\lesssim_f  n^{-(1-p)(\kappa - 1)}   + n^{\frac12 -p + \frac{q}{2} }  + n^{1-2p + q} + o(1).
\end{align*}
Whenever \(p > \frac{q}{2} + \frac12\), the upper bound vanishes as \(n \to \infty\).

Having analyzed all three terms, we conclude that
\eqref{mart_toshow} holds when \(\frac{1}{1+\kappa} < q < \tfrac12\). Such \(q\) exist since \(\kappa > 1\). The proof of Lemma \ref{mart} is complete.
\end{myproof}

The case \(\nu = \mu\) of Theorem \ref{main} now follows from Lemma \ref{mart} exactly as in \cite{dobbs2016}. Below we recall the argument, which uses several facts from stochastic analysis (for proofs of these, see e.g. \cite{rogers2000}). The full result of Theorem \ref{main} follows from the case \(\nu = \mu\) by applying the Portmanteau theorem together with Fact \ref{aimino}.

\subsection{Finishing the proof of Theorem \ref{main}} Assume that \(\nu = \mu\). Since \(\hat{\sigma}_t(f)\) is bounded in \(t\) and independent of \(\chi\), given a standard Brownian motion \((W_t)\), there exists a strong solution to the stochastic differential equation \(d{\chi}(t) = \hat{\sigma}_t(f) \, d W_t\). The solution is unique in law, and we denote its law by \(Q\).

\begin{fact}\label{sol_char} The measure \(Q\) is the unique measure such that  \(Q(\pi_0 = 0) = 1\) and for all \(A \in C_c^{\infty}(\bR) \) the process \((M_t)_{t\in[0,1]}\) is a martingale with respect to \(Q\) and the filtration \((\mathfrak{F}_t)_{t \in [0,1]}\), where \(\mathfrak{F}_t\) is the sigma-algebra on \(C([0,1])\) generated by \(\{ \pi_s \: : \: 0 \le s \le t\} \).
 \end{fact}

By our assumption, \((\bP_n^{\mu})_{n \ge 0} \) is tight, and so there exists a subsequence \((\bP_{n_k}^{\mu})\) which converges in the topology of weak convergence. From Lemma \ref{mart} and Fact \ref{sol_char} we deduce that the weak limit must be \(Q\). This completes the proof of the case \(\nu = \mu\).\\
\indent Then, let \(\nu\) be an arbitrary measure with density \(g = g_1 - g_2\), where \(g_1,g_2 \in \cC_*\), and let \(h\) denote the density of \(\mu\). We denote by \(\bE_n^{\mu,\mu}\) and  \(\bE_n^{\mu,\nu}\) the expectations of \(\bP_n^{\mu,\mu}\) and \(\bP_n^{\mu,\nu}\), respectively. By the Portmanteau theorem and the first part of the proof, it suffices to show that for any bounded and Lipschitz continuous function \(F : C([0,1]) \to \bR\), 
\begin{align}\label{limit_last}
\bE_n^{\mu,\mu}[F] - \bE_n^{\mu,\nu}[F] \to 0.
\end{align}
Since
\begin{align*}
|\bE_n^{\mu,\mu}[F] - \bE_n^{\mu,\nu}[F]| &\le \int |F(\chi_n^{\mu}(x,\cdot)) - F(\chi_n^{\nu}(x,\cdot))| \, d \mu (x)\\
&\le \text{Lip}(F) \sup_{t\in [0,1]} | n^{-\frac12} \mu ( S_n(\cdot, t)) - n^{-\frac12} \nu (S_n(\cdot, t))|,
\end{align*}
\eqref{limit_last} follows by Fact \ref{aimino}:
\begin{align*}
&|n^{-\frac12} \mu ( S_n(\cdot, t)) - n^{-\frac12} \nu (S_n(\cdot, t))| \\
&\le n^{\frac12} \int_0^t |\mu (\bar{f}_{n,\round{ns}}) - m (\bar{f}_{n,\round{ns}})| \, ds + n^{\frac12} \int_0^t |m (\bar{f}_{n,\round{ns}}) - \nu (\bar{f}_{n,\round{ns}})| \, ds \\
&\le \Vert f \Vert_{\infty} n^{\frac12} \left( \int_0^t \Vert \cL_{n,\round{ns}} \cdots \cL_{n,1}  (h -1) \Vert_1 \, ds + \int_0^t \Vert \cL_{n,\round{ns}} \cdots \cL_{n,1}  (g -1) \Vert_1 \, ds \right) \\
&\le \Vert f \Vert_{\infty} n^{-\frac12}  + \sum_{i=1,2} \Vert f \Vert_{\infty} n^{\frac12} \int_0^t \Vert \cL_{n,\round{ns}} \cdots \cL_{n,1}  (g_i - m(g_i)) \Vert_1 \, ds \\
&\lesssim \Vert f \Vert_{\infty} n^{-\frac12} + \Vert f \Vert_{\infty}  (m(g_1) + m(g_2)) n^{-\frac12}.
\end{align*}
We have finished the proof of Theorem \ref{main}. \(\qed\)

\bibliography{qds_pm}{}
\bibliographystyle{plainurl}

\end{document}